\documentclass[11pt,a4paper]{amsart}
\usepackage{amssymb,amsmath,amsthm}

\usepackage[colorlinks=true, pdfstartview=FitV, linkcolor=blue, citecolor=blue, urlcolor=blue,pagebackref=false]{hyperref}

\usepackage{esint}
\usepackage{mathrsfs}

\usepackage{times,latexsym,color}

\newcommand{\BOX}{\ensuremath\Box}

\newtheorem{theorem}{Theorem}
\newtheorem{proposition}{Proposition}
\newtheorem{lemma}[proposition]{Lemma}

\theoremstyle{remark}
\newtheorem{remark}[proposition]{Remark}

\theoremstyle{definition}

\newcommand{\N}{\mathbb{N}}

\newcommand{\C}{\mathbb{C}}
\newcommand{\R}{\mathbb{R}}

\newcommand{\dd}{\,{\rm d}}

\newcommand{\Rel}{\rm Re}

\definecolor{darkgreen}{rgb}{0,0.5,0}
\definecolor{darkblue}{rgb}{0,0,0.7}
\definecolor{darkred}{rgb}{0.9,0.1,0.1}
\definecolor{lightblue}{rgb}{0,0.51,1}

{\vskip\baselineskip\noindent\textbf{Proof of {#1}:}}%
{\hspace*{.1pt}\hspace*{\fill}\BOX\vskip\baselineskip}
{\vskip\baselineskip\noindent\textbf{Proof of Theorem \protect\ref{#1}:}}%
{\hspace*{.1pt}\hspace*{\fill}\BOX\vskip\baselineskip}
{\vskip\baselineskip\noindent\textbf{Proof of Lemma~\protect\ref{#1}:}}%
{\hspace*{.1pt}\hspace*{\fill}\BOX\vskip\baselineskip}
{\vskip\baselineskip\noindent\textbf{Proof of Proposition~\protect\ref{#1}:}}%
{\hspace*{.1pt}\hspace*{\fill}\BOX\vskip\baselineskip}
{\vskip\baselineskip\noindent\textbf{Proof of Theorems \protect\ref{#1} --
\protect\ref{#2}:}}%
{\hspace*{.1pt}\hspace*{\fill}\BOX\vskip\baselineskip}

\begin{document}

\title[Stability of blow-up solutions for Navier-Stokes]{On stability of blow-up solutions of the Burgers vortex type for the Navier-Stokes equations with a linear strain}

\author[Y. Maekawa]{Yasunori Maekawa}
\address[Y. Maekawa]{Kyoto University, Department of Mathematics, Kyoto, Japan}
\email{maekawa@math.kyoto-u.ac.jp}

\author[H. Miura]{Hideyuki Miura}
\address[H. Miura]{Tokyo Institute of Technology, Department of 
Mathematical and Computing Sciences, Tokyo, Japan}
\email{miura@is.titech.ac.jp}

\author[C. Prange]{Christophe Prange}
\address[C. Prange]{Universit\'e de Bordeaux, CNRS, UMR [5251], IMB, Bordeaux, France}
\email{christophe.prange@math.u-bordeaux.fr}

\subjclass[2010]{}
\date{\today}

\maketitle

\noindent {\bf Abstract} We study the three-dimensional Navier-Stokes equations in the presence of the axisymmetric linear strain, where the strain rate depends on time in a specific manner. It is known that the system admits solutions which blow up in finite time and whose profiles are in a backward self-similar form of the familiar Burgers vortices. In this paper it is shown that the existing stability theory of the Burgers vortex leads to the stability of these blow-up solutions as well. The secondary blow-up  is also observed when the strain rate is relatively weak.

\vspace{0.3cm}


\vspace{0.3cm}


\section{Introduction}\label{sec.intro}

One of the important mechanisms in three-dimensional turbulent flows is the vorticity amplification due to stretching and vorticity dissipation due to viscosity.
As a simple model, the vorticity amplification induced by linear straining flows has been widely studied.
A famous example is the Burgers vortex \cite{Bu}, which describes a vortical structure localized in a tubelike domain due to the background axisymmetric linear strain (see also \cite{MKO,GaWa4,Ma1,Ma2} for the study when the linear strain is not necessarily axisymmetric).
Let us consider the three-dimensional Navier-Stokes equations for viscous incompressible flows
\begin{align}\label{eq.3ns}
\partial_t V - \Delta V + \nabla P + V \cdot \nabla V=0\qquad \nabla\cdot V =0, \qquad t>0,\quad x\in \R^3.
\end{align}
The Burgers vortex \cite{Bu} is the steady state solution to \eqref{eq.3ns} of the form
\begin{align}\label{def.burgers.1}
\begin{split}
& V^{rB}(x)= \frac{\gamma}{2} 
\begin{pmatrix}
-x_1\\
-x_2\\
2x_3
\end{pmatrix}
\, + \, \frac{\alpha \big (1- e^{-\frac{\gamma|x'|^2}{4}}\big )}{2\pi |x'|^2}
\begin{pmatrix}
\displaystyle -x_2 \\
\displaystyle x_1\\
0
\end{pmatrix}, \qquad x' = (x_1,x_2),\\
& P^{rB} (x) = -\frac12 |V^{rB} (x)|^2 - \frac{\alpha^2\gamma}{16\pi^2} \int_{\gamma |x'|^2}^\infty \frac{1}{r} (1-e^{-\frac{r}{4}}) e^{-\frac{r}{4}}\, d r,
\end{split}
\end{align}
where $\gamma>0$ and $\alpha\in\R$ are given constants which respectively represent the strain rate and the circulation at infinity. 
The corresponding vorticity field $\Omega^{rB}=\nabla\times V^{rB}$ is given by 
\begin{align}\label{def.burgers.2}
\Omega^{rB}(x)=\begin{split}
\begin{pmatrix}
0\\
0\\
\alpha \gamma g (\gamma^\frac12 x')
\end{pmatrix}, \qquad  
g (X') & = \frac{1}{4\pi} e^{-\frac{|X'|^2}{4}}.
\end{split}
\end{align}
To fix the notation, let 
\begin{equation}\label{e.defU}
U^G(x')=\frac{ 1- e^{-\frac{|x'|^2}{4}}}{2\pi |x'|^2}
\begin{pmatrix}
\displaystyle -x_2 \\
\displaystyle x_1\\
0
\end{pmatrix}\qquad\mbox{and}\qquad G(x')=
\begin{pmatrix}
0\\
0\\
 g (x')
\end{pmatrix}.
\end{equation}
The stability of the Burgers vortex is studied in detail by now, and remarkably, it is stable for any $\gamma>0$ and all $\alpha\in \R$, locally with respect to three-dimensional perturbations and globally with respect to two-dimensional perturbations.  Indeed, the local two-dimensional stability for small circulation is proved by Giga-Kambe \cite{GiKa} and global stability for arbitrary size of circulation is proved by Gallay-Wayne \cite{GaWa1}.
The three-dimensional local stability for small circulation is shown by Gallay-Wayne \cite{GaWa2}, and the restriction on the size of circulation is removed by Gallay-Maekawa \cite{GaMa1} also for the three-dimensional perturbations. The reader is referred to the survey of Gallay-Maekawa \cite{GaMa2} about the existence and the stability problem related to the Burgers vortex.

\smallskip

Although the Burgers vortex presents in a simple way a nontrivial swirling (when $\alpha\ne 0$) flow exhibiting the balance between the vorticity stretching and dissipation, it is a nondecaying (or even growing) solution to the Navier-Stokes equations \eqref{eq.3ns}, and it is well known that we do not have  uniqueness for such flows in general. A typical example is given by the so-called parasitic solutions. Let $\rho:[0,\infty)\rightarrow\R$ with $\rho(0)=\rho_0$ be a bounded differentiable function. Then 
$$
V(x,t)=C\rho(t),\qquad P(x,t)=-\rho'(t)C\cdot x, \qquad C\in \R^3
$$
is a solution to \eqref{eq.3ns} with initial data $C\rho_0$. If one considers the initial data in the form of the linear strain $\rho_0 (-x_1,-x_2,2x_3)^\top$ one finds that
\begin{equation}\label{e.pvlin}
V^{lin}(x,t)=\rho(t) \begin{pmatrix} -x_1\\-x_2\\2x_3 \end{pmatrix},\qquad P(x,t)=-\tfrac12|V^{lin}(x,t)|^2-\tfrac{1}2 \partial_t V^{lin} (x,t)\cdot x
\end{equation}
is a solution to \eqref{eq.3ns}, with a pressure growing quadratically.
In these examples, the function $\rho (t)$ is in principle taken arbitrary, and in particular, one may take it in a singular way so that it blows up at a finite time $T_*>0$.
In view of the linear strain $V^{lin}$ above, it is then natural to look for a solution of the Burgers vortex type but with a time-dependent linear strain which blows up in a finite time, though the function $\rho(t)$ for the strain rate cannot be arbitrary any longer and should be chosen suitably in this case.
Significant contributions in this direction have been made by Moffatt \cite{Mo00} and Ohkitani-Okamoto \cite{OhO}. Indeed, Moffatt \cite{Mo00} provided a family of blowing-up solutions as follows: for $\mu>1$ and $\alpha\in \R$,
\begin{align}
\begin{split}
V^{sB}(x,t)=\ &\frac{\mu}{2(T^*-t)}\begin{pmatrix}
-x_1\\
-x_2\\
2x_3
\end{pmatrix}
\, +\,
\alpha\sqrt{\beta_\mu(t)}U^G\big(\sqrt{\beta_\mu(t)} x'\big),\qquad \beta_\mu (t) = \frac{\mu-1}{T^*-t},\label{e.vsb}\\
P^{sB} (x,t) =\ & -\tfrac{1}2 \partial_t V^{lin} (x,t)\cdot x -\frac12 |V^{sB} (x,t)|^2 - \frac{\alpha^2\beta_\mu (t)}{16\pi^2} \int_{\beta_\mu (t) |x'|^2}^\infty \frac{1}{r} (1-e^{-\frac{r}{4}}) e^{-\frac{r}{4}}\, d r,
\end{split}
\end{align}
where $V^{lin}$ is \eqref{e.pvlin} with $\rho (t) = \frac{\mu}{2(T^*-t)}$ and the formula of the pressure is obtained by using the general identity $u\cdot \nabla u = \frac12 \nabla |u|^2 - u\times \omega$ with $\omega = \nabla \times u$, by also noting the fact that $U^G \times (\nabla \times U^G)=0$ and $\nabla \times V^{lin}=0$.
Indeed, by noticing that the corresponding vorticity field is 
\begin{align}
\Omega^{sB}(x,t)=\ &\alpha\beta_\mu(t)G\big(\sqrt{\beta_\mu(t)}x'\big) \label{e.osb}
\end{align}
with $G$ defined in \eqref{e.defU}, we have
\begin{align*}
&\partial_tV^{sB}+V^{sB}\cdot\nabla V^{sB}-\Delta V^{sB} \\
=\ &  \partial_t V^{lin} +\partial_t \Big (  \alpha\sqrt{\beta_\mu(t)}U^G\big(\sqrt{\beta_\mu(t)} x'\big) \Big )  + \frac12 \nabla |V^{sB}|^2- V^{sB} \times \Omega^{sB} + \nabla \times \Omega^{sB} \\
=\ & \partial_t V^{lin}  + \frac12 \nabla |V^{sB}|^2 -  \alpha\sqrt{\beta_\mu(t)}U^G\big(\sqrt{\beta_\mu(t)} x'\big)  \times \Omega^{sB} \\
\qquad + & \partial_t \Big (  \alpha\sqrt{\beta_\mu(t)}U^G\big(\sqrt{\beta_\mu(t)} x'\big) \Big ) - V^{lin} \times \Omega^{sB} + \nabla \times \Omega^{sB}.
\end{align*}
Then the first three terms in the right-hand side are written as the potential form and hence define the pressure as stated above.
On the other hand, we see from $\partial_t \sqrt{\beta_\mu(t)} = \frac{1}{2(\mu-1)} \beta_\mu(t)^\frac32$, 
\begin{align*}
& \partial_t \Big (  \alpha\sqrt{\beta_\mu(t)}U^G\big(\sqrt{\beta_\mu(t)} x'\big) \Big )  \\
= \ & \frac{\alpha\beta_\mu (t)^\frac32}{\mu-1} \Big ( \frac12 U^G\big(\sqrt{\beta_\mu(t)} x'\big) +\frac12\sqrt{\beta_\mu(t)}x'\cdot\nabla'U^G\big(\sqrt{\beta_\mu(t)}x'\big) \Big )
\end{align*}
and 
\begin{align*}
- V^{lin} \times \Omega^{sB} + \nabla \times \Omega^{sB} & = \frac{\alpha\beta_\mu (t)^2}{2(\mu-1)} \, g (\sqrt{\beta_\mu (t)} x')  \, \begin{pmatrix}
x_2\\
-x_1\\
0
\end{pmatrix} \\
& = \frac{\alpha \beta_\mu (t)^\frac32}{\mu-1} \Delta U^G \big(\sqrt{\beta_\mu(t)}x'\big).
\end{align*}
Thus the conclusion holds from the identity $\Delta U^G (\xi') + \frac12 \xi'\cdot \nabla' U^G (\xi') + \frac12 U^G (\xi') =0$. 
 
The flow \eqref{e.vsb} blows up in a backward self-similar form and has a very similar vortical structure to the Burgers vortex. In this paper we call $V^{sB}$ the singular Burgers vortex. 
Let us notice that, due to the presence of the linear strain, the singular Burgers vortex is out of the theory of Caffarelli-Kohn-Nirenberg \cite{CKN} for $\varepsilon$-regularity for Navier-Stokes equations, and of Ne\v cas-Ru\v zi\v cka-\v Sver\'ak \cite{NRS} and Tsai \cite{T98} for the non existence of backward self-similar solutions.

What is interesting in \eqref{e.vsb} is the restriction $\mu>1$, that is, the strain rate has to be strong enough to define the singular Burgers vortex as a {\it real} flow. For the specific choice of the circulation number $\alpha$, Ohkitani-Okamoto \cite{OhO} provided an alternative interpretation for this time-dependent strain rate in the singular Burgers vortex. Indeed, when $\alpha=\alpha_\mu = \frac{4\pi \mu}{\mu-1}$ one finds that the identity $\frac{\mu}{2(T_*-t)} = \frac{\| \nabla \times V^{sB} (\cdot, t) \|_{L^\infty}}{2}$ holds, that is, the strain rate behaves as if it depends on the unknown variable, i.e., the $L^\infty$ norm of the vorticity field. This gives an interesting perspective in view of the Taylor expansion about $x$ of the velocity around the origin. The analysis in this direction has been developed further by Nakamura-Okamoto-Yagisita \cite{NaOYa} and Okamoto \cite{O}.

From now on we focus our attention on the singular Burgers vortex $V^{sB}$. 
Without loss of generality, we normalize the blow-up time as $1$, i.e., $T^*=1$.
The aim of this paper is to study the asymptotic stability of the explicit blowing-up solutions \eqref{e.vsb}-\eqref{e.osb}. To simplify the notations we set
\begin{align*}
U_\mu(x',t)=\ &\sqrt{\beta_\mu(t)}\, U^G\big(\sqrt{\beta_\mu(t)}x'\big),\\
G_\mu (x',t) =\ &\beta_\mu (t) \, G \big( \sqrt{\beta_\mu (t)} \, x'\big).
\end{align*}
Let us go back to \eqref{eq.3ns} and recall that the vorticity field $\Omega = \nabla \times V$ satisfies 
the equations
\begin{align}\label{eq.3vor}
\partial_t \Omega - \Delta \Omega + V \cdot \nabla \Omega - \Omega \cdot \nabla V =0, \qquad \nabla\cdot\Omega =0,
\end{align}
which are formally equivalent to \eqref{eq.3ns}.
We note that the divergence free condition $\nabla\cdot\Omega=0$ is preserved under the evolution equation in \eqref{eq.3vor}, and thus, if the initial vorticity is divergence free then the second equation in \eqref{eq.3vor} is automatically satisfied. Hence we will always drop the divergence free condition for the vorticity field from now on.
To study the stability of \eqref{e.vsb}-\eqref{e.osb}, we consider the solution $(\Omega, V)$ to \eqref{eq.3vor} of the form $(\Omega, V) = (\Omega^{sB} + \omega, V^{sB} + u)$.
Then  we have from \eqref{eq.3vor} the evolution equations for the perturbation vorticity $\omega$, which reads
\begin{align}\label{eq.perturb.vor}
\begin{split}
\partial_t \omega - \Delta  \omega  + \frac{\mu}{1-t} \Big ( M x \cdot \nabla \omega - M \omega \Big )  + \alpha{\bf \Lambda}_\mu (t) \omega & =- u \cdot \nabla \omega + \omega \cdot \nabla u\\
\omega |_{t=0} & = \omega_0,
\end{split}
\end{align}
and 
\begin{align}
{\bf \Lambda}_\mu (t) \omega & = U_\mu (t) \cdot \nabla \omega + u \cdot \nabla G_\mu (t) - \omega \cdot \nabla U_\mu (t)   - G_\mu (t) \cdot \nabla u.
\end{align}
Here the matrix $M$ is given by 
\begin{align*}
M=
\begin{pmatrix}
-\frac12 & 0 & 0\\
0 & -\frac12 & 0\\
0 & 0 & 1
\end{pmatrix},
\end{align*}
and the perturbation velocity $u$ is formally recovered from the vorticity $\omega$ by the Biot-Savart law 
\begin{align}
u (x,t) =  -\frac{1}{4\pi} \int_{\R^3} \frac{(x-y)\times  \omega (y,t) }{|x-y|^3}  \dd y =  \big ( K_{3D} * \omega (t) \big )  (x),
\end{align}
by assuming a suitable spatial decay on $\omega$.
The goal  of this paper is to study the behavior of $\omega$ up to $T^*=1$ for a suitable class of the initial data $\omega_0$.
In particular, to show the asymptotic stability of the blowing-up solution $(\Omega^{sB}, V^{sB})$, 
we aim at establishing $(1-t)^{\frac12} \| u (t) \|_{L^\infty} \rightarrow 0$ as $t\uparrow 1$.

\smallskip

As in the work of Lundgren \cite{Lu} for the strained Navier-Stokes equations and of Giga-Kohn \cite{GK85} for the nonlinear heat equation, it is convenient to introduce the self-similar variables:
\begin{align*}
\tau = -(\mu-1)\log (1-t),\qquad \xi = \sqrt{\beta_\mu(t)} x, \qquad \beta_\mu (t) = \frac{\mu-1}{1-t}\,,
\end{align*}
and 
\begin{align}\label{def.v.w}
u(x,t) = \sqrt{\beta_\mu (t)}  \mathcal{V} (\xi,\tau) , \qquad \omega (x,t) =  \beta_\mu (t) \mathcal{W}(\xi,\tau).
\end{align}
Then the system for $\mathcal{W}$ is written as  
\begin{align}\label{eq.w}
\begin{split}
\partial_\tau \mathcal{W} - (L_\mu - \alpha {\bf \Lambda} ) \mathcal{W}  & = - \mathcal{V}
 \cdot \nabla \mathcal{W} + \mathcal{W}\cdot \nabla \mathcal{V}  \quad \qquad \tau>0, ~~\xi\in \R^3 ,\\
\mathcal{W}|_{\tau=0}& = \mathcal{W}_0,
\end{split}
\end{align}
where 
\begin{align}\label{def.L_mu}
L_\mu \mathcal{W} & =  \Delta \mathcal{W} + \frac{\xi'}{2}\cdot \nabla' \mathcal{W} - \frac{2\mu+1}{2(\mu-1)} \xi_3 \partial_3 \mathcal{W}  + \frac{1}{\mu-1} ( \mu M-I) \mathcal{W}
\end{align}
and 
\begin{align}\label{def.Lambda}
{\bf \Lambda} \mathcal{W} = U^G \cdot \nabla \mathcal{W} + \mathcal{V} \cdot \nabla G - \mathcal{W} \cdot \nabla U^G - G \cdot \nabla \mathcal{V}.
\end{align}
Here $\mathcal{V}=K_{3D}* \mathcal{W}$,  and $U^G$ and $G$ are defined in \eqref{e.defU}. Hence, by using the self-similar variables $(\xi,\tau)$ the study of the behavior of $\omega$ around the blow-up time $T^*=1$ is translated into the large time behavior of $\mathcal{W}$, and the problem shares common features in essence with the stability problem of the Burgers vortex for which a detailed analysis was done. 
Indeed, if we set the differential operator $L$ as 
\begin{align}\label{def.L}
L w = \Delta w + \frac{\xi'}{2} \cdot \nabla' w - \xi_3 \partial_3 w  + M w,
\end{align}
which is formally obtained by taking the limit $\mu\rightarrow \infty$ in \eqref{def.L_mu},
then the linear operator $L-\alpha {\bf \Lambda}$ is exactly the linearized operator around the Burgers vortex with circulation $\alpha$. It is studied in \cite{GaMa1} for all circulation number $\alpha\in \R$, and the uniform spectral gap of $-\frac12$ is achieved for all $\alpha$ in a suitable functional setting.
Since the only difference between $L_\mu$ and $L$ is the stretching term related to $\xi_3\partial_3$,
it is natural to expect that the operator $L_\mu-\alpha {\bf \Lambda}$ also has a uniform spectral gap under the similar functional framework as in \cite{GaMa1}. This implies that the original blowing-up solution \eqref{e.vsb}-\eqref{e.osb} is asymptotically stable under the small perturbations as long as $\mu>1$. Our main theorem is the following stability result stated in Theorem \ref{theo.main} below, which is in the spirit of stability results for blowing-up solutions of the nonlinear heat equation \cite{MZ97} or of dispersive equations \cite{MM02}. As far as the authors know, this is the first result of this type for the Navier-Stokes equations, though the key stability mechanism is brought by the singular linear strain and hence the spatial growth of the solution plays an important role. 
The function spaces used in the theorem are defined in the paragraph {\bf Notations} below.
It is shown that $L_\mu - \alpha {\bf \Lambda}$ generates a semigroup $e^{\tau (L_\mu -\alpha {\bf \Lambda})}$ in  the function space $\mathbb{X}(m)$ defined in {\bf Notations} and then the mild solution $\mathcal{W}$ to \eqref{eq.w} (with the initial data $\mathcal{W}_0$) is defined as the solution to the integral equation 
\begin{align*}
\mathcal{W} (\tau) & = e^{\tau (L_\mu -\alpha {\bf \Lambda})} \mathcal{W}_0 + \int_0^\tau e^{(\tau-s) (L_\mu -\alpha {\bf \Lambda})} \Big ( - \mathcal{V} \cdot \nabla \mathcal{W} + \mathcal{W}\cdot \nabla \mathcal{V} \Big ) \, d s,\\
\mathcal{V} & =K_{3D}*\mathcal{W}.
\end{align*}
Then we call $\omega$ the mild solution to \eqref{eq.perturb.vor} when $\omega$ is given by the transformation \eqref{def.v.w} for $\mathcal{W}$. 

\begin{theorem}[stability of the singular Burgers vortex]\label{theo.main}
Let $\alpha\in\R$, $m\in (2,\infty]$, $\mu\in(1,\infty)$. Let $\omega_{0,2d}\in L^2_0 (m)$. Then there exists $\varepsilon=\varepsilon(\alpha,m,\mu, \omega_{0,2d})\in(0,1)$ such that the following statement holds. For all divergence-free $\omega_0\in\mathbb X(m)$ satisfying 
$$
\omega_0=(0,0,\omega_{0,2d}(x'))+\omega_{0,3d}(x',x_3),\qquad \|\omega_{0,3d}\|_{\mathbb X(m)}\leq \varepsilon,
$$
there exists a unique mild solution $\omega\in L^\infty_{loc}([0,1);\mathbb X(m))$ to \eqref{eq.perturb.vor} such that $t^\frac12 \partial_x^\beta \omega \in L^\infty_{loc}([0,1);\mathbb X(m))$ with $|\beta|\leq 1$. This solution satisfies
\begin{align*}
\limsup_{t\rightarrow 1} \, (1-t)^{\frac12-\frac{\mu-1}{2}}\|u(\cdot,t)\|_{L^\infty(\R^3)}<\infty,
\end{align*}
where $u=K_{3d}\ast\omega$. Moreover, if $\alpha\neq 0$, letting 
$$
\lambda_i=-\sqrt{\mu-1}\int_{\R^2}x_i\omega_{0,2d}(x')dx',\quad i=\{1,2\},
$$
we have the following asymptotic estimate for $u$:
\begin{equation}\label{e.asy}
\limsup_{t\rightarrow 1}\, (1-t)^{\frac12-\frac{\mu-1}{2}}\Big\|u(\cdot,t)-(1-t)^{\frac{\mu-1}{2}} \sum_{i=1,2}\lambda_i\sqrt{\frac{\mu-1}{1-t}} (\partial_i U^G ) \big(\sqrt{\frac{\mu-1}{1-t}}\cdot\big)\Big\|_{L^\infty(\R^3)} \leq C \varepsilon.
\end{equation}
Here $C$ depends on $\alpha$, $m$, $\mu$, and $\omega_{0,2d}$.
\end{theorem}

\begin{remark}{\rm
(1) Let us notice that $V=V^{sB}+u$ is a solution to the Navier-Stokes system \eqref{eq.3ns}. A way to rephrase the asymptotic expansion \eqref{e.asy} is as follows: 
\begin{align*}
V (x,t) & =V^{sB}(x,t)+\sqrt{\mu-1}(1-t)^{\frac{\mu}{2}-1}\sum_{i=1,2}\lambda_i (\partial_i U^G)\big(\sqrt{\frac{\mu-1}{1-t}}x \big) \\
& \qquad  + O_{L^\infty}\big( \epsilon (1-t)^{\frac{\mu}{2}-1} \big) + o_{L^\infty}\big( (1-t)^{\frac{\mu}{2}-1}\big).
\end{align*}
Hence in the vicinity of $V^{sB}$, we have constructed a whole family of blowing-up solutions. Moreover,  when $\alpha, \lambda_i \ne 0$,  if the strain is sufficiently weak so that $\mu\in(1,2)$, then we identify a leading part (in $\varepsilon$) of the secondary blow-up profile in the sense that the term $\sqrt{\mu-1}(1-t)^{\frac{\mu}{2}-1}\sum_{i=1,2} \lambda_i (\partial_i U^G) \big(\sqrt{\frac{\mu-1}{1-t}}x \big)$ blows up as well in the $L^\infty$ norm.
In fact, as seen in the proof (and the last statement of Theorem \ref{theo.nlstab} and its remark below) the secondary blow up profile is still a linear combination of $\sqrt{\mu-1}(1-t)^{\frac{\mu}{2}-1} (\partial_i U^G)\big(\sqrt{\frac{\mu-1}{1-t}}x \big)$ but in a weaker topology. Precisely, we can show that there exist $d_i\in \R$, $i=1,2$, such that 
\begin{align}\label{e.rem.asy}
\begin{split}
& \lim_{t\rightarrow 1} \, (1-t)^{\frac12-\frac{\mu-1}{2}}\big\|u(\cdot,t) \\
& \qquad -(1-t)^{\frac{\mu-1}{2}} \sum_{i=1,2}(\lambda_i+d_i) \sqrt{\frac{\mu-1}{1-t}} (\partial_i U^G ) \big(\sqrt{\frac{\mu-1}{1-t}}\cdot\big)\big\|_{L^\infty_{x_3}([-N,N]; L^\infty_{x'})} =0,
\end{split}
\end{align}
for any $N>0$. Here $d_i$ satisfies $|d_i|\leq C \|\omega_{0,3d} \|_{\mathbb{X}(m)}\ll 1$ with $C$ depending only on $\alpha$, $m$, $\mu$, and $\omega_{0,2d}$. 

\noindent (2) In Theorem \ref{theo.main} the initial data $\omega_{0,2d}\in L^2_0 (m)$ is taken arbitrary, while we do not have a quantitative information between $\| \omega_{0,2d} \|_{L^2(m)}$ and the small constant $\epsilon$ for the three-dimensional perturbation. On the other hand, for small initial data in $\mathbb{X}(m)$ the condition in Theorem \ref{theo.main} can be stated in a more quantitative way as follows; there exists $\epsilon = \epsilon (\alpha,m,\mu)>0$ such that, if the  (divergence-free) initial data $\omega_0\in \mathbb{X}(m)$ satisfies $\| \omega_0 \|_{\mathbb{X} (m)} \leq \epsilon$, then the stability estimate such as \eqref{e.asy} or \eqref{e.rem.asy} for the solution is verified with a constant $C$ depending only on $\alpha$, $m$, and $\mu$ (here $C$ is taken independently of $\omega_0$, since $\|\omega_0\|_{\mathbb{X}(m)}\leq \epsilon$ and $\epsilon$ is small enough).

\noindent (3) It should be emphasized that the uniqueness of the solution in Theorem \ref{theo.main} is claimed for the equation \eqref{eq.perturb.vor}, and not for the original Navier-Stokes equation \eqref{eq.3ns}. Indeed, as already explained, we do not have the uniqueness of solutions to \eqref{eq.3ns} for nondecaying initial data. Roughly speaking, the uniqueness holds {\it for the class of solutions having the form $V=V^{sB} + u$ with decaying (in the horizontal direction) $u$.} }
\end{remark}

\subsection*{Outline of the paper}

In Section \ref{sec.lin} we handle the stability analysis of the linear equation 
$$
\partial_\tau w - (L_\mu - \alpha {\bf \Lambda} ) w=0.
$$
The goal of the analysis is to extend the arguments of Gallay-Maekawa \cite{GaMa1} to the case of the operator $L_\mu - \alpha {\bf \Lambda}$. 
Section \ref{sec.nl} is devoted to the proof of nonlinear stability. More precisely, we prove that the zero solution of the nonlinear equation \eqref{eq.w} is stable under arbitrarily large two-dimensional perturbations, and small genuinely three-dimensional perturbations, see Theorem \ref{theo.nlstab}. Such a result is in the spirit of \cite{PRST}, though the technique we use is different. Instead of continuing the solution via a blow-up criteria as in \cite{PRST}, we construct a mild solution iteratively until the source becomes small enough for global in time solutions to exist. We conclude this by investigating the existence of a secondary blow-up profile. Theorem \ref{theo.main} immediately follows from the results of Section \ref{sec.nl} by scaling back to the original variables $(x,t)$.

\subsection*{Notations}

As is usual, we always decompose $x=(x',x_3)\in\R^3$ or $\xi=(\xi',\xi_3)\in\R^3$ into horizontal component $x'=(x_1,x_2)\in\R^2$ or $\xi'=(\xi_1,\xi_2)\in\R^2$, and vertical component $x_3$ or $\xi_3$. 
Similary, we write $\nabla'=(\partial_{x_1}, \partial_{x_2})$ and $\Delta'=\partial_{x_1}^2+\partial_{x_2}^2$.
Throughout the paper, we work in the following functional setting. For $m\in [0,\infty]$ set 
\begin{align*}
\rho_m (r) =
\left\{\begin{array}{ll}
 1,   & m=0,\\
 (1+\frac{r}{4m})^m, & 0<m<\infty,\\
 e^{\frac{r}{4}}, & m=\infty.
\end{array}\right.
\end{align*}
For $p\in[1,\infty)$, we define the weighted spaces
\begin{align*}
L^p(m) & = \Big\{w\in L^p (\R^2)~|~\|w\|_{L^p (m)}^2 = \int_{\R^2} |w (x') |^p \rho_m (|x'|^2)^\frac p2  \dd x' \Big\},\\
L^p_0 (m) & = \Big\{ w\in L^p (m)~|~\int_{\R^2} w \dd x' =0\Big\} \quad \mbox{for}\quad m>2-\tfrac2p.
\end{align*}
Moreover, 
we use the following product spaces: 
for $p\in[1,\infty)$ and $m>2-\tfrac2p$ 
\begin{align*}
 X^p (m) =\ & BC (\R; L^p (m)) \,, \qquad X_0^p (m) = BC (\R; L^p_0 (m)),\\
\mbox{with}\quad \|\phi \|_{X^p (m)} =\ & \sup_{x_3\in \R} \| \phi (\cdot,x_3 ) \|_{L^p (m)},\\
\mathbb{X}^p (m) =\ & X^p (m) \times X^p(m) \times X_0^p (m), 
\end{align*}
where ``$BC(\R;Y)$'' denotes 
the space of all bounded and continuous functions from $\R$ to a Banach space $Y$. We denote by $X^p_{loc}(m)$ the subspace of $X^p(m)$ endowed with the topology given by the seminorms $(\|\cdot\|_{X^p_n(m)})_{n\in\N}$ defined by
\begin{equation*}
\|\phi\|_{X^p_n(m)}=\sup_{|x_3|\leq n}\|\phi(\cdot,x_3)\|_{L^p(m)},\quad\forall\phi\in X^p(m),\quad n\in\N.
\end{equation*}
We then set in analogy with above $\mathbb{X}^p_{loc} (m) = X^p_{loc} (m) \times X^p_{loc}(m) \times X_{loc,0}^p (m)$, where $X_{loc,0}(m)$ is $X^p_0(m)$ equipped with the topology of $X^p_{loc}(m)$. When $p=2$ and $m>1$, we use the abbreviations $X(m)$, $X_0(m)$ and $\mathbb X(m)$ and analogously for the ``loc'' versions. These spaces are used in \cite{GaMa1}.

\section{Analysis of the linearized operator}
\label{sec.lin}

This section is centered on the analysis of the semigroup $e^{\tau(L_\mu-\alpha{\bf \Lambda})}$ for the linear evolution. The main result is the following.

\begin{proposition}[linear stability]\label{prop.linstab}
Let $\alpha\in\R$, $m\in (1,\infty]$, $\mu\in(1,\infty)$ and $p\in[1,2]$. Then 
for all $\eta\in(0,\tfrac12]$ with $\eta<\tfrac{m-1}2$, $\kappa \in (0, \eta + \frac{2\mu+1}{2(\mu-1)})$, and $\beta\in\N^3$, there exists a constant $C(\alpha,m,\mu,\kappa,\eta,\beta)<\infty$ such that
\begin{align}
\|\partial_\xi^\beta(e^{\tau(L_\mu-\alpha{\bf \Lambda})}w_0)'\|_{X(m)^2}\leq\ &\frac{Ce^{-(\kappa +\frac{2\mu+1}{2(\mu-1)}\beta_3)\tau}}{a(\tau)^{\frac1p-\frac12+\frac{|\beta|}{2}}}\|w_0\|_{\mathbb X^p(m)}\label{e.linesta}\\
\|\partial_\xi^\beta(e^{\tau(L_\mu-\alpha{\bf \Lambda})}w_0)_3\|_{X(m)}\leq\ &\frac{Ce^{-(\eta+\frac{2\mu+1}{2(\mu-1)}\beta_3)\tau}}{a(\tau)^{\frac1p-\frac12+\frac{|\beta|}{2}}}\|w_0\|_{\mathbb X^p(m)}\label{e.linestb}
\end{align}
for all divergence-free $w_0\in\mathbb X^p(m)$ and $\tau \in(0,\infty)$, where $a(\tau)=1-e^{-\tau}$. Moreover, 
\begin{equation*}
\int_{\R^2}(e^{\tau(L_\mu-\alpha{\bf \Lambda})}w_0)_3(\xi',\xi_3)d\xi'=0,
\end{equation*}
and $\nabla\cdot w_0=0$ implies $\nabla\cdot e^{\tau(L_\mu-\alpha{\bf \Lambda})}w_0=0$ for all $\tau\in(0,\infty)$. If $\eta<\frac12$ then $\kappa$ is taken as $\kappa = \eta + \frac{2\mu+1}{2(\mu-1)}$.
\end{proposition}

Corresponding estimates for the Burgers votex are obatined and used in \cite{GaMa1}. 
The gain in the decay for $\xi_3$ derivatives in \eqref{e.linesta} and \eqref{e.linestb} is due to the commutation property stated in \eqref{e.commx3deriv} below, which is already 
used in \cite{Cr, RoDi, SchRo, GaMa1} . 
This property is an effect of the stretching in the vertical direction due to the structure of the linear strain, which is a key stabilizing effect of the Burgers vortex. Another remark concerns the transient growth. There is a factor $a(\tau)^{-(\frac1p-\frac12+\frac{|\beta|}{2})}$ related to the parabolic-type smoothing effect of $e^{\tau(L_\mu-\alpha{\bf \Lambda})}$, which is large in short time. This factor does not depend on $\alpha$. The constant $C$ in \eqref{e.linesta} and \eqref{e.linestb}, though, gets large when $\alpha\rightarrow\infty$, $\mu\rightarrow 1^+$.
We note that the estimate of $e^{\tau(L_\mu-\alpha{\bf \Lambda})}$ for local time is not difficult to show:
\begin{align}\label{est.local}
\|\partial_\xi^\beta e^{\tau(L_\mu-\alpha{\bf \Lambda})}w_0 \|_{\mathbb{X}(m)}\leq \frac{C}{a(\tau)^{\frac1p-\frac12+\frac{|\beta|}{2}}}\|w_0\|_{\mathbb X^p(m)}, \quad 0<\tau\leq 1, \quad p\in [1,2].
\end{align}
Here $C$ depends only on $\alpha$, $\mu$ and $m$; see, e.g., the argument of \cite[Proposition 4.2]{GaMa1}. Thus, by recalling the semigroup property, we may focus on the estimates \eqref{e.linesta}-\eqref{e.linestb} but only for $\tau\geq 1$ and $p=2$. 

The argument to achieve the linear stability is rather parallel to the one of Gallay-Maekawa in \cite{GaMa1}. The idea is to decompose the full operator $L_\mu-\alpha{\bf \Lambda}$ into a dominant two-dimensional (but vectorial) part and a three-dimensional part whose contribution to the solution decays fast and is negligible in the longtime. Hence we first focus on the operator $L_\mu$. Second, we address the vectorial $2$d problem, i.e. on the action of $L_\mu-\alpha{\bf \Lambda}$ on fields $w=w(x')$ independent of $x_3$. Finally, we analyze the full operator $L_\mu-\alpha{\bf \Lambda}$ using the stretching in the vertical direction which makes the three-dimensional part of the solution decay fast.

\subsection{Analysis of $L_\mu$}

Let us start from the analysis of $L_\mu$. We first rewrite the operator similarly to \cite{GaMa1} so as to make the comparison easier. We first expand \eqref{def.L_mu} as
\begin{equation*}
L_\mu=\Delta'+\frac{\xi'}{2}\cdot\nabla'+\partial_3^2-\frac{2\mu+1}{2(\mu-1)}\xi_3\partial_3+
\left(\begin{array}{ccc}
-\frac{\mu+2}{2(\mu-1)} & 0 & 0\\
0 & -\frac{\mu+2}{2(\mu-1)} & 0\\
0 & 0 & 1
\end{array}\right).
\end{equation*}
Hence, we can write the action of $L_\mu$ in the following form
\begin{equation}
L_\mu w=\left(\begin{array}{c}
L_{\mu,h}w'\\
L_{\mu,3}w_3
\end{array}\right)
=\left(\begin{array}{c}
\mathcal L_hw'+\mathcal L_{3} w'-\big ( 1 + \frac{\mu+2}{2(\mu-1)} \big ) w' \\
\mathcal L_hw_3+\mathcal L_{3}w_3
\end{array}\right),
\end{equation}
where 
\begin{align*}
\mathcal L_h=\ &\Delta'+\frac{\xi'}{2}\cdot\nabla'+1,\\
\mathcal L_{3}=\ &\partial_3^2-\frac{2\mu+1}{2(\mu-1)}\xi_3\partial_3.
\end{align*}
According to \cite[Appendix A]{GaWa2}, the operator $\mathcal L_h$ is the generator of a strongly continuous semigroup in $L^2(m)$ given by the explicit formula
\begin{equation}\label{e.etmathcalLh}
(e^{\tau\mathcal L_h}f)(\xi')=\frac{e^\tau}{4\pi a(\tau)}\int_{\R^2}e^{-\frac{|\xi'-\eta'|^2}{4a(\tau)}}f(\eta'e^{\frac \tau2})d\eta', \qquad a(\tau)=1-e^{-\tau},
\end{equation}
for all $\tau\in(0,\infty)$. Moreover, it is well known that $-\mathcal L_h$ is self-adjoint in $L^2 (\infty)$ and satisfies the following lower bounds (cf. \cite[Lemma 4.7]{GaWa1}):
\begin{align}\label{lower.bound.L}
-\mathcal L_h \geq 0 \quad {\rm in}~L^2(\infty), \qquad -\mathcal L_h \geq \frac12 \quad {\rm in}~L^2_0 (\infty), \qquad -\mathcal L_h \geq 1 \quad {\rm in}~L^2_1 (\infty),
\end{align}
where $L^2_1 (\infty)=\{f\in L^2_0 (\infty)~|~ \int_{\R^2} \xi_j  f  d\xi=0, ~j=1,2\}$.
As for $\mathcal L_3$, it is the generator of a semigroup of contractions in $BC(\R)$ with the following explicit formula, which is derived from \cite[Appendix A]{GaWa2}) for the semigroup generated by the differential operator of the form $\partial_3^2- b \xi_3\partial_3-b$.
\begin{equation}\label{e.etmathcalL3}
(e^{\tau\mathcal L_3}f)(\xi_3)=\big ( \frac{2\chi}{4 \pi (e^{2\chi\tau}-1 )} \big )^\frac12 \int_{\R} \exp  \Big ( -\frac{2\chi |\xi_3-\eta_3|^2}{4 (e^{2\chi \tau} -1)} \Big ) \, f(\eta_3e^{-\chi \tau})d\eta_3,
\end{equation}
for all $\tau\in(0,\infty)$, where 
\begin{equation}\label{e.defgamma}
\chi=\chi_\mu = \frac{2\mu+1}{2(\mu-1)}. 
\end{equation}
Notice that $\chi>1$.  We can rewrite formula \eqref{e.etmathcalL3}:
\begin{equation}\label{e.etmathcalL3'}
(e^{\tau\mathcal L_3}f)(\xi_3)=\left(\frac{\chi}{2\pi a(2\chi \tau)}\right)^\frac12\int_{\R}e^{-\frac{\chi |e^{-\chi \tau}\xi_3-\eta_3|^2}{2a(2\chi \tau)}}f(\eta_3)d\eta_3, \quad a(\tau) = 1-e^{-\tau},
\end{equation}
for all $\tau\in (0,\infty)$. From the previous formula, we immediately obtain the fast temporal decay of derivatives of $e^{\tau\mathcal L_3}f$. Indeed, for all $k\in\N$, there exists a constant $C(k)<\infty$ such that for all $\tau\in(0,\infty)$,
\begin{equation}\label{e.fastdecayL3}
\left\|\partial_3^ke^{\tau\mathcal L_3}f\right\|_{L^\infty(\R)}\leq C\frac{\chi^\frac k2e^{-\chi k\tau}}{a(2\chi \tau)^\frac k2} \| f\|_{L^\infty (\R)},
\end{equation}
where $\chi$ is defined in \eqref{e.defgamma}. Moreover, we also have from \eqref{e.etmathcalL3'},
\begin{align}\label{e.fastdecayL4}
\lim_{\tau \rightarrow \infty} \Big\| e^{\tau\mathcal L_3}f - \left(\frac{\chi}{2\pi}\right)^\frac12\int_{\R}e^{-\frac{\chi |\eta_3|^2}{2}}f(\eta_3)d\eta_3 \Big\|_{L^\infty ([-e^{(\chi-\delta) \tau}, e^{(\chi-\delta))\tau}])}=0,
\end{align}
for any $\delta>0$.

\smallskip

The fast decay \eqref{e.fastdecayL3}, due to the strong stretching in the vertical direction, plays a key role in the fact that the linear evolution becomes independent of $x_3$ at the main order. Hence, the leading dynamics in the longtime is driven by the $2$d vectorial problem, which we analyze in the next subsection.  

\subsection{Localization in the horizontal direction: the $2$d vectorial problem}

In this subsection we analyze the $2$d vectorial problem, which corresponds with the action of $e^{\tau(L_\mu-\alpha{\bf \Lambda})}$ on the functions of the form $w=w(\xi')\in\R^3$, i.e., $\partial_3w=0$. We define for $\alpha\in\R$ and $\mu\in(1,\infty)$,
\begin{equation}\label{e.2dvec}
\mathscr L_{\mu,\alpha} w=\left(\begin{array}{c}
\mathscr L_{\mu,\alpha,h}w'\\
\mathscr L_{\mu,\alpha,3}w_3
\end{array}\right)=\left(\begin{array}{c}
(\mathcal L_h-1-\frac{\mu+2}{2(\mu-1)})w'-\alpha({\bf \Lambda}_1-\tilde{{\bf \Lambda}}_2)w'\\
\mathcal L_hw_3-\alpha({\bf \Lambda}_1+\tilde{{\bf \Lambda}}_3)w_3
\end{array}\right).
\end{equation}
Here
\begin{align*}
{\bf \Lambda}_1w=\ &U^G\cdot\nabla w={(U^G)}'\cdot\nabla' w,\\
\tilde{{\bf \Lambda}}_2w'=\ &w'\cdot\nabla'{(U^G)}',\\
\tilde{{\bf \Lambda}}_3w_3=\ &(K_{2d}\ast w_3)\cdot\nabla' g.
\end{align*}
The operators ${\bf \Lambda}_1$ and $\tilde{{\bf \Lambda}}_3$ come from the vorticity transport, while $\tilde{{\bf \Lambda}}_2$ originates from the vortex stretching term. 
One can check that $L_\mu -\alpha {\bf \bf \Lambda} w = \mathscr L_{\mu,\alpha} w$ when $\partial_3 w=0$; see identity \eqref{e.2d3d} and inequality \eqref{est.3d.perp}. 
Notice that the horizontal $w'$ and the vertical components $w_3$ are completely decoupled. 
This makes the $2$d vectorial problem more tractable than the full  original one. 
The operator $\mathscr L_{\mu,\alpha,3}$ and the associated semigroup in $L^2_0(m)$ are already analyzed in \cite[Section 4]{GaWa1}.
The main result of this section is stated as follows.
\begin{proposition}[$2$d vectorial problem]\label{prop.2dvect}
Let $\alpha\in\R$, $m\in (1,\infty]$, $\mu\in(1,\infty)$. Then for all $\kappa\in(0,1+\tfrac{\mu+2}{2(\mu-1)})$, for all $\eta\in(0,\tfrac12]$ such that $\eta<\tfrac{m-1}2$, there exists a constant $C(\alpha,m,\mu,\kappa,\eta)<\infty$ such that 
\begin{align*}
\|e^{\tau\mathscr L_{\mu,\alpha,h}}w_{0,h}\|_{L^2(m)^2}\leq\ &Ce^{-\kappa \tau}\|w_{0,h}\|_{L^2(m)^2},\\
\|e^{\tau\mathscr L_{\mu,\alpha,3}}w_{0,3}\|_{L^2(m)}\leq\ &Ce^{-\eta \tau}\|w_{0,3}\|_{L^2(m)},
\end{align*}
for all $\tau\in[0,\infty)$ and for all $w_0\in L^2(m)^2\times L^2_0(m)$.
Moreover, if $m>2$ then for $\eta\in (\frac12, 1]$ with $\eta <\frac{m-1}{2}$, 
\begin{align}\label{ex.prop.2dvect}
\Big\|e^{\tau\mathscr L_{\mu,\alpha,3}}w_{0,3} - e^{-\frac{\tau}{2}} \sum_{i=1,2} \theta_i \partial_i g \Big\|_{L^2(m)}\leq\ &C e^{-\eta \tau}\|w_{0,3}\|_{L^2(m)}.
\end{align}
Here $\theta_i = - \int_{\R^2} \xi_i w_{0,3} (\xi') \, d \xi'$ and $g(\xi')=\frac{1}{4\pi} e^{-\frac{|\xi'|^2}{4}}$.
\end{proposition}

Note that the result for $e^{\tau\mathscr L_{\mu,\alpha,3}}w_{0,3}$ in Proposition \ref{prop.2dvect} is due to \cite[Section 4]{GaWa1}, in particular \cite[Proposition 4.12]{GaWa1}. The key observation there is that for suitably large $m$ the spectrum of $\mathscr L_{\mu,\alpha,3}$ in $L^2(m)$ near the imaginary axis consists of the isolated eigenvalues whose eigenfunctions actually belong to $L^2(\infty)$, and thus, the analysis of the large time behavior of $e^{\tau\mathscr L_{\mu,\alpha,3}}$ in $L^2(m)$ is essentially reduced to the analysis in $L^2(\infty)$, in which $\mathcal L_h$ is self-adjoint and moreover ${\bf \Lambda}_1+\tilde{{\bf \Lambda}}_3$ is skew-symmetric; see \cite[Lemma 4.8]{GaWa1}. Then, the lower bounds in \eqref{lower.bound.L} enable us to conclude the expansion  \eqref{ex.prop.2dvect} in $L^2_0(m)$ for large enough $m$, by also using the fact that the eigenspace of the eigenvalue $-\frac12$ of $\mathscr L_{\mu,\alpha,3}$ in $L^2(\infty)$ is spanned by $\partial_i g$, $i=1,2$. On the other hand, the estimate of $e^{\tau\mathscr L_{\mu,\alpha,h}}$ is obtained in the same manner as in \cite[Proposition 3.1]{GaMa1}, as sketched below for reader's convenience.
The following lemma is the key for the study of $e^{\tau\mathscr L_{\mu,\alpha,h}}$. Let $r_{ess}(A; X)$ be the radius of the essential spectrum of a bounded linear operator $A$ on $X$; see \cite[IV-1.20]{ENa}.
\begin{lemma}\label{lem.key2}
Let $m\in (1,\infty]$. Then the following statements hold.

\noindent {\rm (i)} $r_{ess}(e^{\tau\mathscr L_{\mu,\alpha,h}};L^2(m))=\ \begin{cases}
& e^{-(\frac12+\frac{\mu+2}{2(\mu-1)}+\frac m2)\tau}, \quad m\neq\infty\\
& 0, \quad m=\infty
\end{cases}.$

\noindent {\rm (ii)} If $\lambda\in \C$ with ${\Rel} \lambda \geq -\frac12 - \frac{\mu+2}{2(\mu-1)}-\frac{m}{2}$ is an eigenvalue of $\mathscr L_{\mu,\alpha,h}$ in $L^2(m)^2$ then ${\Rel} \lambda\leq -1 - \frac{\mu+2}{2(\mu-1)}$.
\end{lemma}

The proof of Lemma \ref{lem.key2} (i) is identical to the proof of \cite[Proposition 3.3]{GaMa1}. 
Indeed, we see that the operator $\Delta_\alpha(\tau)=e^{\tau\mathscr L_{\mu,\alpha,h}}-e^{\tau(\mathcal L_h-1-\frac{\mu+2}{2(\mu-1)})}$ is compact in $L^2(m)^2$ for any $m\in(1,\infty]$. 
Hence, Weyl's theorem implies that both semigroups have the same essential spectrum and hence have same essential radii: for all $m\in(1,\infty]$, we have 
$r_{ess}(e^{\tau\mathscr L_{\mu,\alpha,h}};L^2(m))^2=\ r_{ess}(e^{\tau(\mathcal L_h-1-\frac{\mu+2}{2(\mu-1)})};L^2(m)^2)$.
Then the fact 
$$r_{ess}(e^{\tau(\mathcal L_h-1-\frac{\mu+2}{2(\mu-1)})};L^2(m))  =\ 
\begin{cases}
& e^{-(\frac12+\frac{\mu+2}{2(\mu-1)}+\frac m2)\tau}, \quad  m\neq\infty\\
& 0, \quad m=\infty
\end{cases},$$ 
which was proved in \cite[Appendix A]{GaWa2}, yields the statement (i) of Lemma \ref{lem.key2}. 
Next, the proof of Lemma \ref{lem.key2} (ii) is sketched below.
By a standard argument \cite[Proposition 3.4 and Section 6.2]{GaMa1}, every eigenfunction in $L^2(m)^2$ associated to an eigenvalue $\lambda$ with ${\Rel(\lambda)}>-\frac12-\frac{\mu+2}{2(\mu-1)}-\frac m2$ belongs to $L^2(\infty)^2$, i.e. has Gaussian decay. It is therefore enough to study the discrete spectrum of $\mathscr L_{\mu,\alpha,h}$ in $L^2(\infty)^2$, for which the same argument as in \cite[Proposition 3.5]{GaMa1} is applied as follows.
The eigenfunction $w'\in L^2(\infty)^2$ associated to the eigenvalue $\lambda$ satisfies, by its definition,
\begin{equation}\label{e.eqeig}
\lambda w'=\mathcal L_hw'-\left(1+\frac{\mu+2}{2(\mu-1)}\right)w'-\alpha {(U^G)}'\cdot\nabla' w'+\alpha w'\cdot\nabla'{(U^G)}'.
\end{equation}
By the direct computation we also have the equations that are respectively satisfied by $\xi'\cdot w'$ and $\nabla' \cdot w'$:
\begin{align}
& \lambda \xi'\cdot w'=\mathcal L_h(\xi'\cdot w')-2\nabla'\cdot w'-\frac12\xi'\cdot w'-\left(1+\frac{\mu+2}{2(\mu-1)}\right)\xi'\cdot w'-\alpha U^G\cdot\nabla'(\xi'\cdot w'), \label{e.eqeig'}\\
& \lambda \nabla'\cdot w'=\mathcal L_h(\nabla'\cdot w')+\frac12\nabla'\cdot w'-\left(1+\frac{\mu+2}{2(\mu-1)}\right)\nabla'\cdot w'-\alpha U^G\cdot\nabla'(\nabla'\cdot w').\label{e.eqeig''}
\end{align}
The upper bound of ${\Rel}\lambda$ is obtained from these identities \eqref{e.eqeig}, \eqref{e.eqeig'}, and \eqref{e.eqeig''}. 
Indeed, testing the equation \eqref{e.eqeig} against $\overline{w'}$ (complex conjugate of $w'$), we obtain
\begin{multline}\label{e.eig1}
{\Rel}\lambda \|w'\|^2=\langle\mathcal L_hw',w'\rangle-\left(1+\frac{\mu+2}{2(\mu-1)}\right)\|w'\|^2\\
+2\alpha{\Rel}\left(\int_{\R^2}e^{\frac{|\xi'|^2}{4}}(\xi'\cdot w')(\xi'^\perp\cdot\overline{w'})\partial_r(u^g)(|\xi'|^2)d\xi'\right),
\end{multline}
where we used the skew-symmetry of $w'\mapsto U^G\cdot\nabla'w'$ in $L^2(\infty)^2$ equipped with the scalar product $\langle\cdot,\cdot\rangle$ defined by
\begin{equation*}
\langle {w'}^1,{w'}^2\rangle=\int_{\R^2}e^{\frac{|\xi'|^2}{4}}{w'}^1\overline{{w'}^2}d\xi'.
\end{equation*} 
Similarly, we have from \eqref{e.eqeig'} and \eqref{e.eqeig''},
\begin{multline}\label{e.eig2}
{\Rel} \lambda \|\xi'\cdot w'\|^2=\langle\mathcal L_h(\xi'\cdot w'),\xi'\cdot w'\rangle-\left(\frac32+\frac{\mu+2}{2(\mu-1)}\right)\|\xi'\cdot w'\|^2\\
-2{\Rel} \langle \nabla'\cdot w',\xi'\cdot w'\rangle,
\end{multline}
and 
\begin{equation}\label{e.eig3}
{\Rel} \lambda\|\nabla'\cdot w'\|^2=\langle\mathcal L_h(\nabla'\cdot w'),\nabla'\cdot w'\rangle-\left(\frac12+\frac{\mu+2}{2(\mu-1)}\right)\|\nabla'\cdot w'\|^2.
\end{equation}
Now suppose that $\nabla'\cdot w'$ is not identically zero. Then \eqref{e.eig3} and \eqref{lower.bound.L} with the fact $\nabla' \cdot w'\in L^2_0(\infty)$ (i.e., $-\mathcal L_h \geq \frac12$ in $L^2_0(\infty)$) imply ${\Rel} \lambda \leq -1 - \frac{\mu+2}{2(\mu-1)}$. Next suppose that $\nabla'\cdot w'$ is identically zero but that $\xi'\cdot w'$ is not identically zero. In this case \eqref{e.eig2} and \eqref{lower.bound.L} (i.e., $-\mathcal L_h \geq 0$ in $L^2(\infty)$) imply ${\Rel} \lambda  \leq -\frac32 - \frac{\mu+2}{2(\mu-1)}$. Finally, suppose that $\xi'\cdot w'$ is identically zero. Then \eqref{e.eig1} and \eqref{lower.bound.L} give the bound ${\Rel} \lambda \leq -1 -\frac{\mu+2}{2(\mu-1)}$. Therefore, we conclude 
\begin{equation}\label{e.boundeig}
{\Rel} \lambda\leq -1-\frac{\mu+2}{2(\mu-1)}.
\end{equation}
The statement of Lemma \ref{lem.key2} (ii) is proved.

\smallskip

The estimate of $e^{\tau\mathscr L_{\mu,\alpha,h}}$ stated in Proposition \ref{prop.2dvect} follows from Lemma \ref{lem.key2} and the standard theory of $C_0$-semigroup \cite[Corollary IV-2.11]{ENa}, and we conclude that the growth bound of $e^{\tau\mathscr L_{\mu,\alpha,h}}$ in $L^2(m)^2$, $m\in (1,\infty]$, is estimated from above by $-1-\frac{\mu+2}{2(\mu-1)}$. Thus, for any $\kappa\in (0, 1+\frac{\mu+2}{2(\mu-1)})$ we have $\|e^{\tau\mathscr L_{\mu,\alpha,h}}w_{0,h}\|_{L^2(m)^2}\leq\ Ce^{-\kappa \tau}\|w_{0,h}\|_{L^2(m)^2}$, which proves Proposition \ref{prop.2dvect}.

\smallskip

We stress that, in Proposition \ref{prop.2dvect}, $w$ is a function of the horizontal variable $\xi'\in\R^2$ only. Notice that we do not yet have the restriction $\kappa$ appearing in Proposition \ref{prop.linstab}. This additional restriction comes from the analysis of the full three-dimensional problem.

\subsection{The regularizing effect of vertical stretching: full $3$d linear stability problem}

In this subsection we complete the proof of Proposition \ref{prop.linstab}. As is mentioned in the beginning of Section \ref{sec.lin}, we may focus on the case $\tau\geq 1$ and $p=2$. 
Our first goal is to show 
\begin{align}\label{est.prop.3dvect.4}
\| e^{(\tau+1) (L_\mu -\alpha {\bf \Lambda})} w_0  \|_{\mathbb{X}(m)} \leq Ce^{-\eta \tau} \| w_0 \|_{\mathbb{X}(m)},
\end{align}
which in particular proves \eqref{e.linestb} with $\beta=0$ and $p=2$ for the vertical component of $e^{\tau (L_\mu -\alpha {\bf \Lambda})} w_0$. 
To this end we see that, due to the stretching effect in the vertical direction, the longtime dynamics of the semigroup $e^{\tau(L_\mu-\alpha\bf\Lambda )}$ is dominated by the $2$d vectorial problem analyzed in Proposition \ref{prop.2dvect}. The following result shows a simple but important stabilizing effect brought by the linear strain to realize this idea.  
\begin{lemma}\label{lem.key1}
We have the following commutation property: for every $\mu>1$ and $\alpha\in\R$,
\begin{equation*}
[\partial_3,L_\mu-\alpha{\bf\Lambda}]=[\partial_3,\mathcal L_3]=-\frac{2\mu+1}{2(\mu-1)}\partial_3.
\end{equation*}
As a consequence, we have for all $\mu>1$,  $\alpha\in\R$, $k\in\N$, for all $\tau\in(0,\infty)$,
\begin{equation}\label{e.commx3deriv}
\partial_3^ke^{\tau(L_\mu-\alpha\bf\Lambda)}=e^{-\frac{2\mu+1}{2(\mu-1)}k\tau} e^{\tau(L_\mu-\alpha\bf\Lambda)}\partial_3^k.
\end{equation}
\end{lemma} 

This property of $\mathcal L_3$ is due to the stretching in the vertical direction and plays a crucial role in reducing the longtime dynamics to the $2$d vectorial problem studied above. Indeed, since it is not difficult to show the naive bound $\| e^{\tau(L_\mu-\alpha\Lambda)}\|_{\mathbb{X}(m) \rightarrow \mathbb{X}(m)} \leq C_1 e^{C_0 \tau}$ for all $\tau>0$, where $C_0$ and $C_1$ may depend on $\alpha$, $\mu$, and $m$, \eqref{e.commx3deriv} gives the bound 
\begin{align*}
\| \partial_3^k e^{(\tau+1) (L_\mu-\alpha\bf\Lambda)} w_0 \|_{\mathbb{X}(m)}  & \leq C_1 e^{C_0 \tau -\frac{2\mu+1}{2(\mu-1)}k\tau} \| \partial_3^k e^{(L_\mu-\alpha\bf\Lambda)} w_0\|_{\mathbb{X}(m)} \\
& \leq C_{1,k}  e^{C_0 \tau -\frac{2\mu+1}{2(\mu-1)}k\tau} \| w_0 \|_{\mathbb{X}(m)}. \quad ({\rm here}~\eqref{est.local} {\rm ~ is ~used})
\end{align*}  
Hence, if $k_0$ is large enough depending on $C_0$, we have 
\begin{align}\label{est.k_0}
\| \partial_3^{k_0} e^{(\tau+1) (L_\mu-\alpha\bf\Lambda)} w_0 \|_{\mathbb{X}(m)}  & \leq C  e^{-\frac{2\mu+1}{2(\mu-1)}\tau} \| w_0 \|_{\mathbb{X}(m)}, \quad \tau > 0.
\end{align}

To obtain the decay estimate \eqref{est.prop.3dvect.4} for $e^{(\tau+1) (L_\mu-\alpha\bf\Lambda)}$, rather than for $\partial_3^{k_0} e^{(\tau+1) (L_\mu-\alpha\bf\Lambda)}$, we decompose the operator $L_\mu-\alpha\bf\Lambda$ into two-dimensional part and three-dimensional part as follows:
\begin{equation}\label{e.2d3d}
(L_\mu-\alpha{\bf \Lambda})w=\mathscr L_{\mu,\alpha}w+\mathcal L_3w
-\alpha\left(\begin{array}{c}
0\\
0\\
{\bf \Lambda}_3w-\tilde{{\bf \Lambda}}_3w_3
\end{array}\right)
+\alpha{\bf \Lambda}_4w,
\end{equation}
where $\mathscr L_{\mu,\alpha}$ is defined in \eqref{e.2dvec} and
\begin{align*}
{\bf \Lambda}_3w-\tilde{{\bf \Lambda}}_3w_3&=\big((K_{3d}\ast w)'-K_{2d}\ast w_3\big)\cdot\nabla' g,\\
{\bf \Lambda}_4w&=g\partial_3(K_{3d}\ast w).
\end{align*}
The proof of the linear stability for the full three-dimensional problem relies on the decomposition \eqref{e.2d3d} and on the following estimates for the three-dimensional part: for all $m\in(1,\infty]$ and $\sigma\in(0,1)$, there exists $C(m,\sigma)$ such that
\begin{align}\label{est.3d.perp}
\begin{split}
\|{\bf \Lambda}_3w-\tilde{{\bf \Lambda}}_3w_3\|_{X(m)}\leq\ & C\big(\|\partial_3w\|_{\mathbb X(m)}+\|w\|_{\mathbb X(m)}^\sigma\|\partial_3w\|_{\mathbb X(m)}^{1-\sigma}\big),\\
\|{\bf \Lambda_4}w\|_{\mathbb X(m)}\leq\ &C\|\partial_3w\|_{\mathbb X(m)}.
\end{split}
\end{align}
These estimates are exactly given in \cite[Proposition 4.5]{GaMa1}, so we omit the proof of \eqref{est.3d.perp}. The three-dimensional part is then treated as a perturbation of the $2$d vectorial problem. The rest of the analysis leading to Proposition \ref{prop.linstab} is rigorously identical to \cite[Section 4]{GaMa1}; namely, it suffices to solve the integral equation
\begin{align}\label{est.prop.3dvect.1}
\begin{split}
\tilde w (\tau) := e^{(\tau+1) (L_\mu -\alpha {\bf \Lambda})} w_0 & = e^{\tau (\mathscr L_{\mu,\alpha} + \mathcal{L}_3  )} \Big (e^{(L_\mu -\alpha {\bf \Lambda})} w_0 \Big ) \\
& \quad - \alpha \int_0^\tau e^{(\tau-s) (\mathscr L_{\mu,\alpha} + \mathcal{L}_3  )} \Bigg \{ \left(\begin{array}{c}
0\\
0\\
{\bf \Lambda}_3 \tilde w-\tilde{{\bf \Lambda}}_3 \tilde w_3
\end{array}\right)
- {\bf \Lambda}_4 \tilde w \Bigg \} \, d s \,,
\end{split}
\end{align}
with the a priori knowledge of the exponential decay of $\partial_3^{k_0}e^{(\tau+1)(L_\mu-\alpha\bf\Lambda)}$ as in \eqref{est.k_0}. Note that the semigroup $e^{\tau (\mathscr L_{\mu,\alpha} + \mathcal{L}_3)}$ is factorized as $e^{\tau \mathcal{L}_3} \otimes  e^{\tau \mathscr L_{\mu,\alpha}}$,
and hence, the estimate of $e^{\tau (\mathscr L_{\mu,\alpha} + \mathcal{L}_3)}$ is a consequence of  \eqref{e.fastdecayL3} and Proposition \ref{prop.2dvect}.
Then,  by applying also \eqref{est.3d.perp} we have for \eqref{est.prop.3dvect.1},
\begin{align}\label{est.prop.3dvect.2}
\begin{split}
& \| \tilde w(\tau ) \|_{\mathbb{X}(m)} \\
&\quad  \leq C e^{-\eta \tau} \| \tilde w_0 \|_{\mathbb{X}(m)} + C \int_0^\tau e^{-\eta (\tau-s)} \big ( \| \partial_3 \tilde w \|_{\mathbb{X}(m)} + \| \tilde w\|_{\mathbb X(m)}^\sigma\|\partial_3 \tilde w\|_{\mathbb X(m)}^{1-\sigma}\big)  (s) \, ds.
\end{split}
\end{align}
Here $\tilde w_0 =  e^{(L_\mu -\alpha {\bf \Lambda})} w_0$ and $C$ depends only on $\alpha$, $\mu$, $m$, and $\sigma \in (0,1)$.
Then the interpolation inequality $\|\partial_3\tilde w\|_{\mathbb X(m)}\leq C \| \tilde w \|_{\mathbb{X}(m)}^{1-\frac{1}{k_0}} \| \partial_3^{k_0} \tilde w \|_{\mathbb{X}(m)}^{\frac{1}{k_0}}$ yields for any $\epsilon\in (0,1)$, 
\begin{align*}
\| \tilde w(\tau ) \|_{\mathbb{X}(m)} \leq C e^{-\eta \tau} \| \tilde w_0 \|_{\mathbb{X}(m)} + C \int_0^\tau e^{-\eta (\tau-s)} \big (\epsilon \| \tilde w \|_{\mathbb{X}(m)} + C_\epsilon \|\partial_3^{k_0} \tilde w\|_{\mathbb X(m)} \big)  (s) \, ds.
\end{align*}
Note that the term $\int_0^\tau e^{-\eta(\tau-s)} \| \partial_3^{k_0} \tilde w (s) \|_{\mathbb{X}(m)} d s$ is bounded from above by $Ce^{-\eta \tau} \| w_0 \|_{\mathbb{X}(m)}$, in virtue of \eqref{est.k_0} and $\eta\in (0,\frac12]$. Then, by taking $\epsilon$ small enough, one can show that $\| \tilde w(\tau) \|_{\mathbb{X}(m)} \leq C e^{-\frac{\eta}{2} \tau} \| \tilde w_0 \|_{\mathbb{X}(m)} + C e^{-\eta \tau} \| w_0 \|_{\mathbb{X}(m)}$. This estimate combined with the local (in time) estimate implies, in the end, 
\begin{align}\label{est.nonsharp}
\| e^{\tau (L_\mu -\alpha {\bf \Lambda})} w_0 \|_{\mathbb{X}(m)} \leq C e^{-\frac{\eta}{2}\tau} \| w_0 \|_{\mathbb{X}(m)}, \quad \tau>0.
\end{align}
The decay rate is then improved as follows. From \eqref{e.commx3deriv} and \eqref{est.nonsharp} we have 
\begin{align}\label{est.prop.3dvect.3}
\| \partial_3 e^{(\tau+1) (L_\mu -\alpha {\bf \Lambda})} w_0 \|_{\mathbb{X}(m)} \leq C e^{-\frac{\eta}{2}\tau - \frac{2\mu+1}{2(\mu-1)}\tau} \| \partial_3 \tilde w_0 \|_{\mathbb{X}(m)} \leq Ce^{-\frac{\eta}{2}\tau - \frac{2\mu+1}{2(\mu-1)}\tau} \| w_0 \|_{\mathbb{X}(m)}.
\end{align}
Then, \eqref{est.prop.3dvect.2} with $\sigma$ close to $0$ and \eqref{est.prop.3dvect.3} imply \eqref{est.prop.3dvect.4} and hence \eqref{e.linestb} for $\beta=0$.
For $\beta \neq 0$, \eqref{e.linestb} follows from  \eqref{est.local},  
\eqref{est.prop.3dvect.4}, and the identity 
\begin{equation*}
\partial^\beta_{\xi} e^{\tau(L_\mu-\alpha\bf\Lambda)}w_0
=e^{-\frac{2\mu+1}{2(\mu-1)}\beta_3 (\tau-\frac12)} \partial^{\beta'}_{\xi'}
e^{\frac12(L_\mu-\alpha\bf\Lambda)}
e^{(\tau-1)(L_\mu-\alpha\bf\Lambda)}
\partial^{\beta_3}_{3}e^{\frac12(L_\mu-\alpha\bf\Lambda)} w_0.
\end{equation*}

Next, to show \eqref{e.linesta}, it suffices to consider the case $\beta=0$.
Fix a given number $\kappa \in (0,\eta +\tfrac{2\mu+1}{2(\mu-1)})$.
As for the horizontal component of $e^{(\tau+1) (L_\mu -\alpha {\bf \Lambda})} w_0$, denoted by $\tilde w'(\tau)$, we have again from \eqref{est.prop.3dvect.1} and Proposition \ref{prop.2dvect},
\begin{align}\label{est.prop.3dvect.5}
\| \tilde w'(\tau) \|_{X(m)^2} \leq C e^{-\kappa' \tau} \| \tilde w_0' \|_{X(m)^2} + C \int_0^\tau e^{-\kappa' (\tau-s)}  \| \partial_3 \tilde w \|_{\mathbb{X}(m)} (s)  d s.
\end{align}
Here $\kappa' \in (0,1+\tfrac{\mu+2}{2(\mu-1)})$ is taken so that $\kappa'>\kappa$, which is possible since $\eta\in (0,\frac12]$ and $\kappa \in (0,\eta +\tfrac{2\mu+1}{2(\mu-1)})$. Note that \eqref{e.commx3deriv} and \eqref{est.prop.3dvect.4} imply $\| \partial_3  e^{(\tau+1) (L_\mu -\alpha {\bf \Lambda})} w_0  \|_{\mathbb{X}(m)} \leq C e^{-\eta \tau - \frac{2\mu+1}{2(\mu-1)}\tau} \| w_0 \|_{\mathbb{X}(m)}$.
Thus \eqref{est.prop.3dvect.5} gives the bound $\| \tilde w'(\tau) \|_{X(m)^2} \leq C e^{-\kappa \tau}\| w_0 \|_{\mathbb{X}(m)}$, as desired. In the case $\eta<\frac12$ one can take the above $\kappa'$ as $\kappa'> \eta +\tfrac{2\mu+1}{2(\mu-1)}$, which gives the decay rate $e^{-(\eta +\tfrac{2\mu+1}{2(\mu-1)})\tau}$ when $\eta<\frac12$. This concludes the proof of Proposition \ref{prop.linstab}.

\subsection{Long time asymptotics for the full $3$d linearized problem}
 
In this subsection we show the asymptotic estimate of $e^{\tau(L_\mu -\alpha {\bf \Lambda})}$ for large $\tau$. The main result is the estimate \eqref{est.prop.2dvect.5} below. 
Let $m>2$ and let us introduce the projection $\mathcal{P}_1$ as 
\begin{align}\label{est.prop.2dvect.2}
\mathcal{P}_1 f = \sum_{i=1,2} \theta_i [f_3] \partial_i G , \qquad G = \begin{pmatrix} 0 \\ 0 \\ \displaystyle \frac{1}{4\pi} e^{-\frac{|\xi'|^2}{4}}  \end{pmatrix}, \qquad f\in \mathbb{X}(m).
\end{align}
Here $\theta_i [f_3] (\xi_3) = - \int_{\R^2} \xi_i f_3 (\xi', \xi_3) \, d \xi'$.
Then $\mathcal{P}_1$ commutes with $e^{\tau \mathcal L_3}$ and we also note that $\mathcal{P}_1$ is motivated by the eigenprojection of $\mathscr L_{\mu,\alpha,3}$ for the eigenvalue $-\frac12$.
In particular, we have 
\begin{align}\label{est.prop.2dvect.3}
\begin{split}
\Big ( \mathcal{P}_1 e^{\tau (\mathscr L_{\mu,\alpha} + \mathcal L_3)} f  \Big ) (\xi)
& = e^{-\frac{\tau}{2}} \sum_{i=1,2}
\big ( e^{\tau \mathcal{L}_3}  \theta_i [f_3] \big ) (\xi_3) \,  \partial_i G (\xi'),  \qquad f\in \mathbb{X} (m).
\end{split}
\end{align}
Set $w (\tau) = e^{\tau (L_\mu -\alpha {\bf \Lambda})} w_0$, $w_0\in \mathbb{X}(m)$, which satisfies the statement in Proposition \ref{prop.linstab}: since $m>2$ we have 
\begin{align}\label{est.long.lin1}
\| w(\tau) \|_{\mathbb{X}(m)} \leq Ce^{-\frac\tau2} \| w_0 \|_{\mathbb{X}(m)}, \qquad \| \partial_3 w (\tau) \|_{\mathbb{X}(m)} \leq \frac{C e^{-(\frac12+\frac{2\mu+1}{2(\mu-1)})\tau}}{a(\tau)^\frac12} \| w_0 \|_{\mathbb{X}(m)}.
\end{align}  
As in \eqref{est.prop.3dvect.1}, $w$ satisfies the formula
\begin{align}\label{est.long.lin2}
w (\tau) & = e^{\tau (\mathscr L_{\mu,\alpha} + \mathcal{L}_3  )} w_0
- \alpha \int_0^\tau e^{(\tau-s) (\mathscr L_{\mu,\alpha} + \mathcal{L}_3  )} \Big ( \left(\begin{array}{c}
0\\
0\\
{\bf \Lambda}_3w-\tilde{{\bf \Lambda}}_3w_3
\end{array}\right)
- {\bf \Lambda}_4 w \Big ) \, d s \,,
\end{align}
Let $\eta\in (\frac12,1]$ with $\eta<\frac{m-1}{2}$.
Proposition \ref{prop.2dvect} and \eqref{e.fastdecayL3} together with \eqref{est.3d.perp} yield
\begin{align*}
\| (I-\mathcal{P}_1) w(\tau) \|_{\mathbb{X}(m)} & \leq C e^{-\eta \tau} \| w_0 \|_{\mathbb{X}(m)} \\
& \quad + C \int_0^\tau e^{-\eta(\tau-s)} \big ( \| \partial_3 w\|_{\mathbb{X}(m)} + \| w\|_{\mathbb{X}(m)}^\sigma \| \partial_3 w\|_{\mathbb{X}(m)}^{1-\sigma} \big ) (s)  d s.
\end{align*}
Since $\frac{2\mu+1}{2(\mu-1)}>1$, we then have from \eqref{est.long.lin1} and by taking $\sigma$ close to $0$,
\begin{align}
\| (I-\mathcal{P}_1) w(\tau) \|_{\mathbb{X}(m)} & \leq C e^{-\eta \tau} \| w_0 \|_{\mathbb{X}(m)} + C \int_0^\tau e^{-\eta(\tau-s)} \frac{e^{-\frac32s}}{a(s)^\frac12} d s \, \| w_0 \|_{\mathbb{X}(m)} \nonumber \\
& \leq C e^{-\eta \tau} \| w_0 \|_{\mathbb{X}(m)}. \label{est.prop.2dvect.4}
\end{align}
We also observe from \eqref{est.long.lin2} and the definition of $\mathcal{P}_1$,
\begin{align}\label{est.long.lin3}
\mathcal{P}_1 w(\tau) =  e^{-\frac{\tau}{2}}  \sum_{i=1,2}  \Big ( e^{\tau \mathcal{L}_3}  \theta_i [w_{0,3}]  - \alpha \int_0^\tau e^{\frac{s}{2}} e^{(\tau-s) \mathcal{L}_3}  h_i (s) \, d s \Big ) \partial_i G, 
\end{align}
with 
\begin{align*}
h_i (\xi_3,s) = - \int_{\R^2} \xi_i \Big ( {\bf \Lambda}_3w-\tilde{{\bf \Lambda}}_3w_3 - ({\bf \Lambda}_4 w )_3 \Big ) (\xi', \xi_3, s)  d \xi'.
\end{align*}
Notice that $h_i$ satisfies from \eqref{est.3d.perp} that
\begin{align*}
\| h_i (s) \|_{L^\infty (\R)} & \leq C \big ( \| \partial_3 w\|_{\mathbb{X}(m)} + \| w\|_{\mathbb{X}(m)}^\sigma \| \partial_3 w\|_{\mathbb{X}(m)}^{1-\sigma} \big ) (s) \\
& \leq \frac{Ce^{-\frac32 s}}{a(s)^\frac12} \| w_0 \|_{\mathbb{X}(m)},  \qquad ({\rm from ~\eqref{est.long.lin1} ~and}~ \sigma~{\rm is~taken~ as~ close ~to}~0)
\end{align*} 
which implies from \eqref{e.fastdecayL4}, with $\chi=\frac{2\mu+1}{2(\mu-1)}$,
\begin{align}\label{est.long.lin4}
\begin{split}
& \lim_{\tau \rightarrow \infty} \Big \| \int_0^\tau e^{\frac{s}{2}} e^{(\tau-s) \mathcal{L}_3}  h_i (s) \, d s \\
& \quad -\left(\frac{\chi}{2\pi}\right)^\frac12 \int_0^\infty e^{\frac{s}{2}} \int_{\R}e^{-\frac{\chi |\eta_3|^2}{2}} h_i (\eta_3,s) d\eta_3 d s  \Big \|_{L^\infty ([-e^{(\chi-\delta) \tau}, e^{(\chi-\delta)\tau}])}=0,
\end{split}
\end{align}
for any small $\delta>0$.
We also have again from \eqref{e.fastdecayL4} that
\begin{align}\label{est.long.lin5}
\lim_{\tau \rightarrow \infty} \Big \|  e^{\tau \mathcal{L}_3}  \theta_i [w_{0,3}]   - \left(\frac{\chi}{2\pi}\right)^\frac12  \int_{\R}e^{-\frac{\chi |\eta_3|^2}{2}} \theta_i [w_{0,3}] (\eta_3) d\eta_3 \Big \|_{L^\infty ([-e^{(\chi-\delta) \tau}, e^{(\chi-\delta)\tau}])}=0.
\end{align}
Thus, \eqref{est.prop.2dvect.4}, \eqref{est.long.lin3}, \eqref{est.long.lin4}, and \eqref{est.long.lin5} give the following asymptotic estimate:
\begin{align}\label{est.prop.2dvect.5}
\begin{split}
& \lim_{\tau\rightarrow \infty} \Big\| e^{\frac{\tau}{2}} e^{\tau (L_\mu -\alpha {\bf \Lambda})} w_0 - \left(\frac{\chi}{2\pi}\right)^\frac12 \sum_{i=1,2} \Big ( \int_{\R}e^{-\frac{\chi |\eta_3|^2}{2}} \theta_i [w_{0,3}] d\eta_3\\
& \qquad  - \alpha \int_0^\infty e^{\frac{s}{2}} \int_{\R}e^{-\frac{\chi |\eta_3|^2}{2}} h_i (\eta_3,s) d\eta_3 d s \Big ) \partial_i G\Big\|_{L^\infty ([-e^{(\chi-\delta) \tau}, e^{(\chi-\delta)\tau}]; L^2(m)^3)}=0,
\end{split}
\end{align}
for any small $\delta>0$.
The similar convergence is valid also for $\nabla w' (\tau)$.
Note that the coefficient in the expansion in \eqref{est.prop.2dvect.5} satisfy 
\begin{align*}
\left(\frac{\chi}{2\pi}\right)^\frac12  \Big | \int_{\R}e^{-\frac{\chi |\eta_3|^2}{2}} \theta_i [w_{0,3}] d\eta_3 - \alpha \int_0^\infty e^{\frac{s}{2}} \int_{\R}e^{-\frac{\chi |\eta_3|^2}{2}} h_i (\eta_3,s) d\eta_3 d s \Big | \leq C \| w_0 \|_{\mathbb{X}(m)},
\end{align*}
with $C$ depending only on $\alpha$, $m$, and $\mu$.

\section{Nonlinear stability}
\label{sec.nl}

This section is devoted to the analysis of the longtime behavior of the nonlinear system 
\eqref{eq.w}. We stress that the longtime behavior of solutions $w$ to \eqref{eq.w} immediately translates into information about the behavior of perturbations of the singular Burgers vortex near the blow-up time $T^*=1$. Hence Theorem \ref{theo.main} is a direct consequence of the results of this section.

\smallskip

We prove that the solution $0$ of 
\begin{align}\label{eq.wbis}
\partial_\tau \mathcal{W} - (L_\mu - \alpha {\bf \Lambda} ) \mathcal{W}  & = - \mathcal{V} \cdot \nabla \mathcal{W} + \mathcal{W} \cdot \nabla \mathcal{V} \quad\mbox{on}\quad (0,\infty)\times\R^3 
\end{align}
is asymptotically stable with respect to perturbations of the divergence-free initial data of the form
\begin{equation*}
\mathcal{W}_0=W_{0,2d}+W_{0,3d}(\xi',\xi_3),\qquad\mbox{with}\quad W_{0,2d}=\left(\begin{array}{c} 
0\\
0\\
w_{0,2d}(\xi')
\end{array}\right)
\end{equation*}
where $w_{0,2d}$ is a scalar field in $L^2_0(m)$ of arbitrary size and $W_{0,3d}$ is small in $\mathbb X(m)$. The strategy we use is reminiscent of the paper \cite{PRST}. We first study the stability of \eqref{eq.wbis} with respect to arbitrarily large $2$d perturbations $W_{0,2d}$ in $L^2_0(m)$. For this we rely on the result of \cite{GaWa1} about the longtime behavior of Navier-Stokes equations in $\R^2$. Let us call $(V_{2d},W_{2d})$ the solution to 
\begin{align}\label{eq.w2d}
\begin{split}
\partial_\tau W_{2d} - (L_\mu - \alpha {\bf \Lambda} ) W_{2d}  & = - V_{2d} \cdot \nabla W_{2d} + W_{2d}\cdot \nabla V_{2d}  \quad \qquad \tau>0, ~~\xi\in \R^3 ,\\
W_{2d}|_{\tau=0}& = W_{0,2d}.
\end{split}
\end{align}
This is in fact a two-dimensional problem and the solution is of the form $\begin{pmatrix} 0 \\ 0 \\ w_{2d} (\xi',\tau) \end{pmatrix}$, and hence is reduced to the scalar equation for $w_{2d} (\xi',\tau)$ discussed in \cite{GaWa1}. We will state the result for $w_{2d}$ in Subsection \ref{subsec.2d} below.  
Then we study the stability of the solution $0$ to the perturbed system around $(V_{2d}, W_{2d})$ with  the small perturbation $W_{0,3d}$ in $\mathbb X(m)$ of the initial data.
Precisely, we consider the solution $\mathcal{W}=w+W_{2d}$ to \eqref{eq.wbis} with the initial data $\mathcal{W}_0$, and thus, the equation for $w$ reads
\begin{align}\label{eq.wpert}
\begin{split}
&\partial_\tau w - (L_\mu - \alpha {\bf \Lambda} ) w\\
&\quad=  - v \cdot \nabla w + w\cdot \nabla v
 -V_{2d}\cdot\nabla w-v\cdot\nabla W_{2d}+W_{2d}\cdot\nabla v+w\cdot\nabla V_{2d}=:F\\
&\qquad\qquad\qquad\qquad\qquad\qquad\qquad\qquad\qquad \mbox{on}\quad (0,\infty)\times\R^3 
\end{split}
\end{align}
where $v=K_{3D}*w$ is given by the Biot-Savart law. For this we use a fixed point argument 
treating the linear terms
\begin{equation*}
 -V_{2d}\cdot\nabla w-v\cdot\nabla W_{2d}+W_{2d}\cdot\nabla v+w\cdot\nabla V_{2d}
\end{equation*}
perturbatively.

\smallskip

Theorem \ref{theo.main} is a reformulation in the original variables of the following stability theorem.

\begin{theorem}[nonlinear stability]\label{theo.nlstab}
Let $\alpha\in\R$, $m\in (1,\infty]$ and $\mu\in(1,\infty)$. Let  $w_{0,2d}\in L^2_0(m)$.
For all $\eta\in(0,\frac12]$ such that $\eta<\frac{m-1}2$, there exists $\varepsilon_0(\alpha,m,\mu,\eta, w_{0,2d} )\in(0,\infty)$, such that for all $W_{0,3d}\in \mathbb X(m)$, the condition 
\begin{equation*}
\|W_{0,3d}\|_{\mathbb X(m)}\leq \varepsilon_0,
\end{equation*}
implies there exists a unique mild solution $\mathcal{W}\in L^\infty((0,\infty);\mathbb X(m))\cap C^0([0,\infty);\mathbb X_{loc}(m))$ to \eqref{eq.wbis} with the initial data $W_0=\begin{pmatrix} 0 \\ 0 \\ w_{0,2d} \end{pmatrix} + W_{0,3d}$ satisfying 
\begin{equation*}
\big\|\partial^\beta_\xi \big ( \mathcal{W} (\cdot, \tau) -W_{2d}(\cdot,\tau) \big )\big\|_{\mathbb X(m)}\leq \frac{C}{a(\tau)^{\frac{|\beta|}{2}}} e^{-\eta\tau} \| W_{0,3d} \|_{\mathbb{X}(m)},\qquad\forall \tau\in(0,\infty).
\end{equation*}
Here $|\beta|\leq 1$ and $C$ depends on $\alpha$, $m$, $\mu$, $\eta$, and $w_{0,2d}$.
Moreover, if $m>2$ then there exist $d_i\in \R$, $i=1,2$, such that 
\begin{align}\label{est.theo.nlstab.1}
\lim_{\tau \rightarrow \infty} \sum_{|\beta|\leq 1}  \Big\| \partial_\xi ^\beta \big ( e^\frac{\tau}{2}  \mathcal{W} (\cdot, \tau) - \sum_{i=1,2} (\lambda_i + d_i) \partial_i G (\cdot) \big )  \Big\|_{L^\infty_{\xi_3}([-e^{(\chi-\delta)\tau},e^{(\chi-\delta)\tau}]; L^2(m)^3)} =0,
\end{align}
for any $\delta>0$. Here $\lambda_i = -\int_{\R^2} \xi_i w_{0,2d} (\xi') \, d \xi'$ and $d_i$ satisfies $|d_i|\leq C \| W_{0,3d} \|_{\mathbb{X}(m)}$, and $\chi=\frac{2\mu+1}{2(\mu-1)}>1$. If $W_{0,3d}=0$ then $d_i=0$ and \eqref{est.theo.nlstab.1} is valid in $L^\infty_{\xi_3}(\R; L^2(m)^3)$. 
\end{theorem}

\begin{remark}{\rm From the Biot-Savart law and \eqref{est.theo.nlstab.1} we have 
\begin{align}\label{est.rem.theo.nlstab.1}
\lim_{\tau\rightarrow \infty} \Big\| e^\frac{\tau}{2} \mathcal{V} (\cdot,\tau) - \sum_{i=1,2} (\lambda_i + d_i) \partial_i U^G (\cdot)   \Big\|_{L^\infty_{\xi_3}([-e^{(\chi-\delta)\tau},e^{(\chi-\delta)\tau}]; L^\infty_{\xi'})} =0,
\end{align}
for any $\delta>0$. Since $\chi>1$ we obtain \eqref{e.rem.asy} by rescaling back to the original variable.
}
\end{remark}

\subsection{Global stability with respect to $2$d perturbations}\label{subsec.2d}

We look for a solution $W_{2d}$ of \eqref{eq.w2d} in the form
\begin{equation*}
W_{2d}(\xi',\tau)=\left(\begin{array}{c}
0\\
0\\
w_{2d}(\xi',\tau)
\end{array}\right).
\end{equation*}
We notice that
\begin{equation*}
L_\mu W_{2d}=\left(\Delta'
+\tfrac12\xi'\cdot\nabla'+1\right)\left(\begin{array}{c}
0\\
0\\
w_{2d}
\end{array}\right)=\left(\begin{array}{c}
0\\
0\\
L_{2d} w_{2d}
\end{array}\right)
\end{equation*}
and 
\begin{equation*}
{\bf \Lambda} W_{2d}=(U^G)'\cdot\nabla'\left(\begin{array}{c}
0\\
0\\
w_{2d}
\end{array}\right)
+v'\cdot\nabla'\left(\begin{array}{c}
0\\
0\\
g(\xi')
\end{array}\right)=\left(\begin{array}{c}
0\\
0\\
{\bf \Lambda}_{2d} w_{2d}
\end{array}\right)
\end{equation*}
so that \eqref{eq.w2d} becomes the scalar equation
\begin{align}\label{eq.w2dscalar}
\begin{split}
\partial_\tau w_{2d} - (L_{2d} - \alpha {\bf \Lambda}_{2d} ) w_{2d}  & = - v_{2d} \cdot \nabla w_{2d}  \quad \qquad \tau>0, ~~\xi\in \R^3 ,\\
w_{2d}|_{\tau=0}& = w_{0,2d}.
\end{split}
\end{align}
The stability of \eqref{eq.w2dscalar} was studied in full details. We recall here the result of \cite{GaWa1}. Pay attention to the fact that the result in \cite{GaWa1} is for the stability of the Oseen vortex with vorticity $g$, while we focus on the stability of the $0$ solution of \eqref{eq.w2dscalar}. In the end both viewpoints are of course equivalent since $w$ solving \eqref{eq.w2dscalar} is nothing but the perturbation of the Oseen vortex.

\begin{proposition}[{\cite[Proposition 1.5, 1.6 and Proposition 4.14]{GaWa1}}]\label{prop.summaryGaWa}
Let $m\in(1,\infty]$.  
For all divergence-free $w_{0,2d}\in L^2_0(m)$, there exists a unique global solution $w_{2d}\in C^0([0,\infty);L^2(m))$ and
\begin{equation*}
\|w_{2d}(\cdot,\tau)\|_{L^2(m)}\longrightarrow 0\quad\mbox{when}\quad \tau\rightarrow\infty.
\end{equation*}
Moreover, for all $\eta\in(0,\tfrac12]$ such that $\eta<\tfrac{m-1}{2}$, for all $\beta\in\N^2$, $|\beta|\leq 1$, for all $w_{0,2d}\in L^2_0(m)$, there exists a constant $C(\eta,\beta,w_{0,2d})<\infty$ such that 
\begin{equation*}
\|\partial_{\xi'}^\beta w_{2d}(\cdot,\tau)\|_{L^2(m)}\leq \frac{Ce^{-\eta\tau}}{a(\tau)^\frac{|\beta|}{2}},
\end{equation*}
for all $\tau\in (0,\infty)$. 
\end{proposition}

The result stated in \cite{GaWa1} does not explicitly include the case $m=\infty$. However, the global stability for the case $m=\infty$ follows from the global stability for the case $m<\infty$. Indeed, the perturbed velocity $\| v_{2d} (\cdot,\tau) \|_{\infty}$ becomes small in long time in virtue of the result for $m<\infty$, which is enough to close the estimate of $w_{2d}$ in $L^2(\infty)$. This argument is explicitely written in \cite[Section 5.1]{Ga18}. 
Furthermore, let us remark that the fact that $w_{0,2d}$ is of integral zero over $\R^2$, i.e. $w_{0,2d}\in L^2_0(m)$ is crucial. Otherwise, the perturbation $w$ would not go to zero, since $\alpha g+w_{2d}$ would go to the Oseen vortex with total circulation $\alpha+\int_{\R^2}w_{0,2d}(\xi')d\xi'$. Notice that for $m>2$, we have the decay
\begin{equation*}
\|w_{2d}(\cdot,\tau)\|_{L^2(m)}\lesssim e^{-\frac\tau 2},
\end{equation*}
which is sharp if one of the following moments is different from zero
\begin{equation*}
\int_{\R^2}\xi_1w_{0,2d}(\xi')d\xi'\neq 0\quad\mbox{or}\quad \int_{\R^2}\xi_2w_{0,2d}(\xi')d\xi'\neq 0.
\end{equation*}
The result of Gallay and Wayne goes even beyond the statement in Proposition \ref{prop.summaryGaWa}, in the sense that the next term in the asymptotics of $w_{2d}$ is given, see \cite[display (80)]{GaWa1} and Subsection \ref{sec.sec} below.

\subsection{Local stability with respect to $3$d perturbations}

Our next aim is to analyze the stability of \eqref{eq.wpert} for general, but small in $\mathbb X(m)$, initial perturbations $W_{0,3d}$. The driving terms of the evolution are the linear terms $L_\mu-\alpha{\bf \Lambda}$. Hence, we treat the remaining terms as perturbations. Indeed the nonlinear terms
\begin{equation*}
- v \cdot \nabla w + w\cdot \nabla v
\end{equation*}
are small because the initial data is small, and the linear terms
\begin{equation*}
-V_{2d}\cdot\nabla w-v\cdot\nabla W_{2d}+W_{2d}\cdot\nabla v+w\cdot\nabla V_{2d}
\end{equation*}
are small in longtime thanks to the results of Proposition \ref{prop.summaryGaWa}.

\smallskip

It follows from this observation that the solution $w$ is a solution of the integral equation
\begin{align*}
w(\cdot,\tau)=\ &e^{\tau(L_\mu-\alpha{\bf \Lambda})}W_{0,3d}+\sum_{i=1}^2\int_0^\tau e^{(\tau-\sigma)(L_\mu-\alpha{\bf \Lambda})}N_i(w,w)d\sigma\\
&+\sum_{i=3}^6\int_0^\tau e^{(\tau-\sigma)(L_\mu-\alpha{\bf \Lambda})}N_i(w)d\sigma\\
=\ &e^{\tau(L_\mu-\alpha{\bf \Lambda})}W_{0,3d}+\sum_{i=1}^2\phi_i(w,w)+\sum_{i=3}^6\phi_i(w)\\
=\ &\Phi(w)(\cdot,\tau),
\end{align*}
where
\begin{align*}
N_1(w,\tilde w)=-(K_{3d}*w)\cdot\nabla \tilde w,&\quad N_2(w,\tilde w)=w\cdot\nabla (K_{3d}*\tilde w),\\
N_3(w)=-V_{2d}\cdot\nabla w,&\quad N_4(w)=-(K_{3d}*w)\cdot\nabla W_{2d},\\
N_5(w)=W_{2d}\cdot\nabla (K_{3d}*w),&\quad N_6(w)=w\cdot\nabla V_{2d},
\end{align*}
with $K_{3d}$ the three-dimensional Biot-Savart kernel. Therefore, $w$ is a fixed point of the map $\Phi$. 

\subsubsection*{Function space and a priori bounds} The linear evolution satisfies the estimates of Proposition \ref{prop.linstab}. Let $m\in(1,\infty]$ and $\eta\in(0,\tfrac12]$ such that $\eta<\tfrac{m-1}{2}$. We introduce the Banach space $\mathbb U$ defined as follows
\begin{multline*}
\mathbb U=\Big\{w\in L^\infty((0,\infty);\mathbb X(m))\cap C^0([0,\infty);\mathbb X_{loc}(m)); \nabla\cdot w(\cdot,\tau)=0,\quad \forall \tau\in(0,\infty),\\
\|w\|_{\mathbb U}=\sum_{|\beta|\leq 1}\sup_{\tau>0} a(\tau)^\frac{|\beta|}{2}e^{\eta\tau}\|\partial_\xi^\beta w(\cdot,\tau)\|_{\mathbb X(m)}<\infty
\Big\}.
\end{multline*}
Notice that the definition of $\mathbb U$ depends on $\eta$ and $m$, though this dependence is not explicitly written. To deal with the linear perturbation terms $\phi_i$ for $i\in\{3,\ldots\ 6\}$ for the short time existence, we also introduce the space $\mathbb U_T$ for given $T\in(0,\infty)$ defined as follows
\begin{multline*}
\mathbb U_T=\Big\{w\in L^\infty((0,T);\mathbb X(m))\cap C^0([0,T);\mathbb X_{loc}(m)); \nabla\cdot w(\cdot,\tau)=0,\quad \forall \tau\in(0,T),\\
\|w\|_{\mathbb U_T}=\sum_{|\beta|\leq 1}\sup_{\tau\in(0,T)} a(\tau)^\frac{|\beta|}{2}e^{\eta\tau}\|\partial_\xi^\beta w(\cdot,\tau)\|_{\mathbb X(m)}<\infty
\Big\}.
\end{multline*}
From Proposition \ref{prop.linstab}, we now easily obtain the following estimate for the linear evolution
\begin{equation*}
\|e^{\tau(L_\mu-\alpha{\bf \Lambda})}W_{0,3d}\|_{\mathbb U}\leq C_1\|W_{0,3d}\|_{\mathbb X(m)},
\end{equation*}
with a constant $C_1(m,\eta,\mu)<\infty$.

\smallskip

The nonlinear terms $\phi_1$ and $\phi_2$ can be handled exactly as in \cite[Section 5]{GaMa1} using \cite[Corollary 2.4]{GaMa1}. For all $p\in(1,2)$, there exists a constant $C(m,p)<\infty$ such that for all $w,\ \tilde w\in\mathbb X(m)$, for all $i\in\{1,2\}$,
\begin{equation*}
\|N_i(w,\tilde w)\|_{X^p(m)^3}\leq C\|w\|_{\mathbb X(m)}\|\nabla \tilde w\|_{\mathbb X(m)}.
\end{equation*}
Hence following the estimates in \cite[p. 503]{GaMa1}, for all $p\in(1,2)$, there exists a constant $C(m,p,\eta)<\infty$, for all $\beta\in\N^3$, $|\beta|\leq 1$, for all $i\in\{1,2\}$, for all $w,\ \tilde w\in\mathbb X(m)$,
\begin{equation}\label{e.estphii}
\left\|\partial_x^\beta\phi_i(w,\tilde{w})\right\|_{\mathbb X(m)}\leq \frac{Ce^{-\eta\tau}}{a(\tau)^{\frac1p+\frac{|\beta|}2-1}}\|w\|_{\mathbb U}\|\tilde w\|_{\mathbb U}\leq \frac{Ce^{-\eta\tau}}{a(\tau)^{\frac{|\beta|}2}}\|w\|_{\mathbb U}\|\tilde w\|_{\mathbb U},
\end{equation}
using that $p>1$ and that $a$ is bounded for the last inequality. Moreover, for $w\in\mathbb U$, we have $\nabla\cdot w=0$ by definition, which implies
\begin{equation*}
\int_{\R^2}(N_1(w,w)+N_2(w,w))_3d\xi'=\int_{\R^2}\nabla'\cdot(w'v_3-v'w_3)d\xi'=0.
\end{equation*}
Therefore, $N_1(w,w)+N_2(w,w)\in\mathbb X^p(m)$ and hence 
\begin{equation*}
\int_{\R^2}(\phi_1(w,w)+\phi_2(w,w))_3d\xi'=0.
\end{equation*}
By estimate \eqref{e.estphii}, $\phi_1(w,w)+\phi_2(w,w)$ in addition belongs to $\mathbb U$ and
\begin{equation*}
\left\|\phi_1(w,w)+\phi_2(w,w)\right\|_{\mathbb U}\leq C_2\|w\|_{\mathbb U}^2,
\end{equation*}
with a constant $C_2(m,p,\eta)<\infty$.

\smallskip

We now turn to the linear terms $\phi_i(w)$, for $i\in\{3,\ldots\, 6\}$. These terms are negligible in longtime due to the exponential decay exhibited in Proposition \ref{prop.summaryGaWa}. In short time however, there is no obvious decay. We will use the factor $a(\tau)^{-\frac1p+1}$ neglected in the estimate \eqref{e.estphii} of the nonlinear terms in order to gain smallness. Thanks to Proposition \ref{prop.summaryGaWa}, we can estimate directly: for all $p\in(1,2)$, there exists a constant $C(m,p,\eta,W_{0,2d})<\infty$ such that for all $w\in\mathbb X(m)$, for $i\in\{3, 5\}$,
\begin{equation*}
\|N_i(w)\|_{X^p(m)^3}\leq C\|w_{2d}\|_{L^2(m)}\|\nabla w\|_{\mathbb X(m)}\leq Ce^{-\eta\tau}\|\nabla w\|_{\mathbb X(m)}\leq \frac{Ce^{-2\eta\tau}}{a(\tau)^\frac12}\|w\|_{\mathbb U},
\end{equation*}
and for $i\in\{4,6\}$,
\begin{equation*}
\|N_i(w)\|_{X^p(m)^3}\leq C\|\nabla w_{2d}\|_{L^2(m)}\|w\|_{\mathbb X(m)}\leq \frac{Ce^{-\eta\tau}}{a(\tau)^\frac12}\|w\|_{\mathbb X(m)}\leq \frac{Ce^{-2\eta\tau}}{a(\tau)^\frac12}\| w\|_{\mathbb U}.
\end{equation*}
Therefore, one can estimate $\phi_i$, for $i\in\{3,\ldots\, 6\}$ in the exact same way as the non\-linear terms $\phi_i$ for $i\in\{1,2\}$ above. This yields for all $p\in(1,2)$, there exists a constant $C(m,p,\eta, W_{0,2d})<\infty$, for all $\beta\in\N^3$, $|\beta|\leq 1$, for all $i\in\{3,\ldots\, 6\}$, for all $w\in\mathbb X(m)$,
\begin{equation}\label{e.estphiilin}
\left\|\partial_x^\beta\phi_i(w)\right\|_{\mathbb X(m)}\leq \frac{Ce^{-\eta\tau}}{a(\tau)^{\frac1p+\frac{|\beta|}2-1}}\|w\|_{\mathbb U},
\end{equation}
Contrary to the nonlinear terms, here we keep the factor $a(\tau)^{-\frac1p+1}$ which is used to give smallness in short time, using the immediate inequality
\begin{equation}\label{e.ineqa}
a(\tau)\leq\tau,\qquad\forall\tau\in(0,\infty).
\end{equation}
It remains to see that the horizontal moments of the third components of $N_3(w)+N_6(w)$ on the one hand and $N_4(w)+N_5(w)$ on the other hand are zero. Indeed, \begin{align*}
\int_{\R^2}(N_3(w)+N_6(w))_3d\xi'
=\ &\int_{\R^2}-V_{2d}'\cdot\nabla' w_3+w'\cdot\nabla' V_{2d,3}d\xi'\\
=\ &-\int_{\R^2}\nabla'\cdot(V_{2d}'w_3)d\xi'=0
\end{align*}
and
\begin{align*}
\int_{\R^2}(N_4(w)+N_5(w))_3d\xi'=\ &\int_{\R^2}-v'\cdot\nabla' W_{2d,3}+W_{2d,3}\partial_3v_3d\xi'
\\
=\ &-\int_{\R^2}\nabla\cdot(v'W_{2d,3})d\xi'=0.
\end{align*}
Thus, $\phi_3(w)+\ldots\ \phi_6(w)$ belongs to $\mathbb U$ and we have the estimate
\begin{equation*}
\big\|\sum_{i=3}^6\phi_i(w)\big\|_{\mathbb U}\leq C_3\|w\|_{\mathbb U}
\end{equation*}
as well as
\begin{equation*}
\big\|\sum_{i=3}^6\phi_i(w)\big\|_{\mathbb U_T}\leq C_3T^{1-\frac1p}\|w\|_{\mathbb U_T}
\end{equation*}
using \eqref{e.ineqa}, for all $T\in(0,\infty)$, with a constant $C_3(m,p,\eta, W_{0,2d})<\infty$.

\subsubsection*{Short time existence on $(0,T)$} Let $p\in(1,2)$. To put it in a nutshell, for $T\in (0,\infty)$, we proved for all $w\in\mathbb U$
\begin{align*}
\|\Phi(w)\|_{\mathbb U}\leq\ & C_1\|W_{0,3d}\|_{\mathbb X(m)}+C_2\|w\|_{\mathbb U}^2+C_3\|w\|_{\mathbb U},\\
\|\Phi(w)\|_{\mathbb U_T}\leq\ & C_1\|W_{0,3d}\|_{\mathbb X(m)}+C_2\|w\|_{\mathbb U_T}^2+C_3T^{1-\frac1p}\|w\|_{\mathbb U_T}
\end{align*}
and for all $w,\ \tilde w\in\mathbb U_T$
\begin{equation*}
\|\Phi(w)-\Phi(\tilde w)\|_{\mathbb U_T}\leq C_2(\|w\|_{\mathbb U_T}+\|\tilde w\|_{\mathbb U_T})\|w-\tilde w\|_{\mathbb U_T}+C_3T^{1-\frac1p}\|w-\tilde w\|_{\mathbb U_T}
\end{equation*}
The goal is now to choose $T\in(0,\infty)$, $\delta>0$ and $K>0$ such that if $W_{0,3d}\in B_{\mathbb X(m)}(0,\delta)$, then $\Phi$ maps $B_{\mathbb U_T}(0,K)$ into itself and is a contraction on this ball. In order to realize this, it is enough to take $T$, $\delta$ and $K$ such that
\begin{equation*}
C_3T^{1-\frac1p}<\frac12,\quad 2C_2K<\frac12\quad\mbox{and}\quad \delta<\frac{K}{4C_1}.
\end{equation*}
Therefore, there exists a unique fixed point $w$ for $\Phi$ in $\mathbb U_T$.

\subsubsection*{Longtime existence and decay} It remains to extend the solution on $(0,\infty)$ and to prove decay estimates for $w$. Let $T$ be fixed for the whole discussion as above
\begin{equation*}
C_3T^{1-\frac1p}<\frac12.
\end{equation*}
We use the fact that $W_{2d}$ and $V_{2d}$ decay fast in time. One of the arguments we use repeatedly is the following. For $j\in\N$, $j\geq 1$, we consider $w_j$ defined by $w_j(\cdot,\tau)=w(\cdot,\tau+jT)$ for all $\tau\in(0,\infty)$ and which solves
\begin{align}\label{eq.wperttilde}
\begin{split}
\partial_\tau w_j - (L_\mu - \alpha {\bf \Lambda} ) w_j & =  - v_j \cdot \nabla w_j + w_j\cdot \nabla v_j\\
 -V_{2d}^{(j)}\cdot\nabla w_j-v_j&\cdot\nabla W_{2d}^{(j)}+W_{2d}^{(j)}\cdot\nabla v_j+w_j\cdot\nabla V_{2d}^{(j)},\qquad \tau>0, ~~\xi\in \R^3, \\
w_j|_{\tau=0}& = w(\cdot,jT),
\end{split}
\end{align}
where $v_j=K_{3d}*w_j$, $V_{2d}^{(j)}=V_{2d}(jT+\cdot)$ and $W_{2d}^{(j)}=W_{2d}(jT+\cdot)$. We rewrite this system as an integral equation, defining the shifted version of the mapping $\Phi$, $\Phi^{(j)}$, in an obvious way. Notice that the smallness of $V_{2d}^{(j)}$ and $W_{2d}^{(j)}$ follows easily from the exponential decay of $V_{2d}$ and $W_{2d}$. Hence we obtain, following the same estimates as above for $\Phi$ but with $V_{2d}$ and $W_{2d}$ replaced by $V_{2d}^{(j)}$ and $W_{2d}^{(j)}$, for all $w\in\mathbb U$
\begin{align*}
\|\Phi^{(j)}(w)\|_{\mathbb U}\leq\ & C_1\|w(\cdot,jT)\|_{\mathbb X(m)}+C_2\|w\|_{\mathbb U}^2+C_3e^{-\eta jT}\|w\|_{\mathbb U}
\end{align*}
and for all $w,\ \tilde w\in\mathbb U$
\begin{equation*}
\|\Phi^{(j)}(w)-\Phi^{(j)}(\tilde w)\|_{\mathbb U}\leq C_2(\|w\|_{\mathbb U}+\|\tilde w\|_{\mathbb U})\|w-\tilde w\|_{\mathbb U}+C_3e^{-\eta jT}\|w-\tilde w\|_{\mathbb U}.
\end{equation*}
For all $w,\ \tilde w\in\mathbb U_T$,
\begin{align}\label{e.estwj}
\begin{split}
\|\Phi^{(j)}(w)\|_{\mathbb U_T}\leq\ & C_1\|w(\cdot,jT)\|_{\mathbb X(m)}+C_2\|w\|_{\mathbb U_T}^2+C_3e^{-\eta jT}T^{1-\frac1p}\|w\|_{\mathbb U_T},\\
\|\Phi^{(j)}(w)-\Phi^{(j)}(\tilde w)\|_{\mathbb U_T}\leq\ & C_2(\|w\|_{\mathbb U_T}+\|\tilde w\|_{\mathbb U_T})\|w-\tilde w\|_{\mathbb U_T}+C_3e^{-\eta jT}T^{1-\frac1p}\|w-\tilde w\|_{\mathbb U_T}\\
\leq\ &C_2(\|w\|_{\mathbb U_T}+\|\tilde w\|_{\mathbb U_T})\|w-\tilde w\|_{\mathbb U_T}+C_3T^{1-\frac1p}\|w-\tilde w\|_{\mathbb U_T}
\end{split}
\end{align}
The idea is to iterate the short time construction as long as $C_3e^{-\eta jT}\geq\frac12$. Let $k$ be the least integer such as 
\begin{equation}\label{e.condk}
C_3e^{-\eta kT}<\frac12.
\end{equation}
Notice that $k(m,\eta,p)<\infty$. Let $K_{\infty}$ be such that
\begin{equation}\label{e.condKk+1}
2C_2K_{\infty}<\frac12.
\end{equation}
We subsequently iterate $k$ times the short time construction in order to have a solution in $(0,kT)$. Then, we can construct the solution directly in $(kT,\infty)$ based on the fact that the perturbative terms are small uniformly in that time interval. For all $j\in\{0,k+1\}$, let 
\begin{equation*}
K_j=\frac{1}{4C_1}K_{j+1}\qquad\mbox{so that}\quad K_{j}=\big(\frac{1}{4C_1}\big)^{k+2-j}K_{\infty}.
\end{equation*}
Let us explain the induction. First for $k=0$, we construct the solution in $(0,T)$. This was above, in the local in time construction of a solution. We take $\delta=K_0$ and $K=K_1$. By our choice of parameters and estimates \eqref{e.estwj}, we obtain the existence of a unique solution to \eqref{eq.wpert}, which satisfies
\begin{equation*}
\|w\|_{\mathbb U_T}\leq K_1\quad\mbox{and moreover, by definition of}\ \mathbb U_T,\quad\|w(\cdot,\tau)\|_{\mathbb X(m)}\leq K_1 e^{-\eta\tau},\ \forall\tau\in(0,T).
\end{equation*}
We can now iterate the construction. Let $j\in\{1,k-1\}$. Assume that we have a solution $w$ on $(0,jT)$ such that
\begin{equation}\label{e.reca}
w\in L^\infty((0,jT);\mathbb X(m))\cap C^0([0,jT);\mathbb X_{loc}(m));\ \nabla\cdot w(\cdot,\tau)=0\ \forall \tau\in(0,T)
\end{equation}
and 
\begin{equation}\label{e.recb}
\|w(\cdot,\tau)\|_{\mathbb X(m)}\leq K_j e^{-\eta\tau},\ \forall\tau\in(0,T).
\end{equation}
Then, we aim at extending the solution on $(0,(j+1)T)$. In order to do so, we consider the shifted solution $w_j$ to \eqref{eq.wperttilde}. By \eqref{e.estwj} we have
\begin{align*}
\|\Phi^{(j)}(w)\|_{\mathbb U_T}\leq\ & C_1K_je^{-\eta jT}+C_2K_{j+1}\|w\|_{\mathbb U_T}+C_3e^{-\eta jT}T^{1-\frac1p}K_{j+1},\\
\leq\ & \frac{3e^{-\eta jT}}{4}K_{j+1}+\frac{1}{4}\|w\|_{\mathbb U_T}\\
\leq\ & K_{j+1}\\
\|\Phi^{(j)}(w)-\Phi^{(j)}(\tilde w)\|_{\mathbb U_T}\leq\ & (2C_2K_{j+1}+C_3T^{1-\frac1p})\|w-\tilde w\|_{\mathbb U_T}\\
\leq\ & A\|w-\tilde w\|_{\mathbb U_T},
\end{align*}
with $A<1$. Therefore, there exists a unique fixed point $w_j$ such that
\begin{equation*}
\|w_j\|_{\mathbb U_T}\leq K_{j+1}.
\end{equation*}
Furthermore, thanks to the bound
\begin{equation*}
\|w_j\|_{\mathbb U_T}=\|\Phi^{(j)}(w)\|_{\mathbb U_T}\leq \frac{3e^{-\eta jT}}{4}K_{j+1}+\frac{1}{4}\|w\|_{\mathbb U_T}
\end{equation*}
we obtain
\begin{equation*}
\|w_j\|_{\mathbb U_T}\leq K_{j+1}e^{-\eta jT},
\end{equation*}
which implies 
\begin{equation*}
\|w(\cdot,\tau)\|_{\mathbb X(m)}\leq K_{j+1}e^{-\eta\tau},\quad \forall\tau\in(0,(j+1)T),
\end{equation*}
where $w$ is the concatenation of $w$ defined on $(0,jT)$ and $w_j(\cdot,\cdot-jT)$ on $(jT,(j+1)T)$. Hence, we have the recurrence hypothesis \eqref{e.reca} and \eqref{e.recb} at rank $j+1$.

\smallskip

It remains to construct a solution on $(kT,\infty)$. This is can now be done in one step since we do not need the smallness of the parameter $T$ any longer to make $C_3e^{-\eta kT}K_{\infty}$ small. By the condition \eqref{e.condk}, we have
\begin{equation*}
C_3e^{-\eta kT}K_{\infty}<\frac12 K_{\infty}.
\end{equation*}
Therefore, we can easily show that there exists a unique $w_{k+1}$ on $(0,\infty)$ such that
\begin{equation*}
\|w_{k+1}\|_{\mathbb U}\leq K_{\infty}.
\end{equation*}
Concatenating this solution with the one on $(0,kT)$ and using the previous estimates, we arrive at the existence of $w$ such that
\begin{equation*}
w\in L^\infty\big((0,\infty);\mathbb X(m)\big)\cap C^0\big([0,\infty);\mathbb X_{loc}(m)\big),\quad \nabla\cdot w(\cdot,\tau)=0,\ \forall \tau\in(0,\infty)
\end{equation*}
and
\begin{equation*}
\|w(\cdot,\tau)\|_{\mathbb X(m)}\leq K_\infty e^{-\eta\tau},\quad\forall\tau\in(0,\infty).
\end{equation*}
This ends the proof of Theorem \ref{theo.nlstab}. 

\begin{remark}
Notice that the size of $3$d perturbations $\|W_{0,3d}\|_{\mathbb X(m)}$ that are allowed in the argument above depends on $w_{0,2d}$ through the constant $C_3$. Indeed, 
\begin{equation*}
\|W_{0,3d}\|_{\mathbb X(m)}\leq\delta\simeq \frac1{4C_2}e^{-(k+2)\log(4C_1)}.
\end{equation*}
From our choice of parameters, we have the following rough estimates: $T\simeq (2C_3)^{-\frac1{1-\frac1p}}$, hence from \eqref{e.condk}
\begin{equation*}
k\simeq \frac1{\eta T}\log(2C_3)\simeq \frac{1}{\eta}(2C_3)^\frac1{1-\frac1p}\log(2C_3).
\end{equation*}
\end{remark}

\subsection{Secondary blow-up profile}
\label{sec.sec}

Let $\alpha\in\R$, $\alpha\neq 0$ be fixed. Let $\nu\in (\frac12,1]$ and $m>2\nu +1$.
Note that we have constructed the solution to \eqref{eq.wbis} of the form 
\begin{align*}
\mathcal{W} (\xi,\tau) = \begin{pmatrix} 0 \\ 0 \\ w_{2d} (\xi',\tau)\end{pmatrix} + w (\xi,\tau), 
\end{align*}
with 
\begin{align}\label{proof.sec.sec.1}
\begin{split}
\sup_{\tau\geq 0} \sum_{|\beta|\leq 1} e^{\frac{\tau}{2}} a(\tau)^\frac{|\beta|}{2} \| \partial_\xi^\beta w_{2d} (\tau) \|_{L^2 (m)} & \leq C (W_{0,2d}), 
\\
\sup_{\tau\geq 0} \sum_{|\beta|\leq 1} e^{\frac{\tau}{2}} a(\tau)^\frac{|\beta|}{2} \| \partial_\xi^\beta w(\tau ) \|_{\mathbb{X}(m)} & \leq C \| W_{0,3d} \|_{\mathbb{X}(m)} \ll 1.
\end{split}
\end{align}
Notice in addition that by Gallay-Wayne \cite[Eq.(80)]{GaWa1} with the smoothing effect of the system we have 
\begin{align}\label{proof.sec.sec.2}
\begin{split}
& \sup_{\tau\geq 0} \sum_{|\beta|\leq 1} e^{\nu \tau} a(\tau)^\frac{|\beta|}{2} \Big\| \partial_\xi^\beta \Big ( w_{2d} (\tau) - e^{-\frac{\tau}{2}} \sum_{i=1,2} \lambda_i \partial_i g  \Big )\Big\|_{L^2 (m)} \leq C (W_{0,2d}), \\
& \lambda_i = - \int_{\R^2} \xi_i w_{0,2d} (\xi') \, d \xi'.
\end{split}
\end{align}
Then \eqref{proof.sec.sec.1} and \eqref{proof.sec.sec.2} yield from the Biot-Savart law,
\begin{align}\label{proof.sec.sec.0}
\limsup_{\tau\rightarrow \infty} \, \Big\| e^\frac{\tau}{2}  \mathcal{V} (\tau) - \sum_{i=1,2} \lambda_i \partial_i U^G \Big\|_{L^\infty} \leq C \| W_{0,3d} \|_{\mathbb{X}(m)}.
\end{align}
Theorem \ref{theo.main} now follows from rescaling Theorem \ref{theo.nlstab} and the above estimate \eqref{proof.sec.sec.0} for the velocity.

It remains to show the last statement of Theorem \ref{theo.nlstab}.
We observe from \eqref{eq.wpert} that $w$ is the solution to 
\begin{align}\label{proof.sec.sec.3}
w(\tau) = e^{\tau (L_\mu - \alpha {\bf \Lambda})} W_{0,3d} + \int_0^\tau e^{(\tau-s) (L_\mu - \alpha {\bf \Lambda})} F(s) \, d s\,,
\end{align}
where we already know that $F(s)$ satisfies the estimate 
\begin{align}\label{proof.sec.sec.4}
\| F(s) \|_{\mathbb{X}(m)} \leq \frac{C}{a (s)^\frac34} e^{-\frac34s}\| W_{0,3d} \|_{\mathbb{X}(m)}, \qquad C=C(W_{0,2d})\,.
\end{align}
Then \eqref{est.prop.2dvect.5} implies that for each $f\in \mathbb{X}(m)$ with $m>2$ there exist $c_i[f]\in \R$, $i=1,2$, such that 
\begin{align}\label{proof.sec.sec.5}
\lim_{\tau\rightarrow \infty} \sum_{|\beta|\leq 1} \Big\| \partial_\xi^\beta \Big ( e^{\frac{\tau}{2}} e^{\tau (L_\mu - \alpha {\bf \Lambda})} f - \sum_{i=1,2} c_i[f]  \partial_i G \Big ) \Big\|_{L^\infty_{\xi_3} ([-e^{(\chi-\delta)\tau},e^{(\chi-\delta)\tau}]; L^2 (m)^3)}=0,
\end{align}
for any $\delta>0$. Here $\chi=\frac{2\mu+1}{2(\mu-1)} >1$ and the coefficients $c_i[f]$ satisfy
$$|c_i[f]|\leq C \| f\|_{\mathbb{X}(m)}. $$
Thus, by combining \eqref{proof.sec.sec.4}, \eqref{proof.sec.sec.5}, and the estimate $\| \partial^\beta_\xi e^{\tau (L_\mu - \alpha {\bf \Lambda})}  f\|_{\mathbb{X}(m)} \leq C a(\tau)^{-\frac{|\beta|}{2}} e^{-\frac{\tau}{2}} \| f\|_{\mathbb{X}(m)}$, we have from the integral equation \eqref{proof.sec.sec.3} that 
\begin{align}\label{proof.sec.sec.6}
\begin{split}
& \lim_{\tau\rightarrow \infty} \sum_{|\beta \leq 1}\Big\| \partial_\xi^\beta \Big ( e^{\frac{\tau}{2}} w(\tau) \\
& \qquad  - \sum_{i=1,2} \Big ( c_i [W_{0,3d}] + \int_0^\infty e^\frac{s}{2} c_i[F(s)] \, d s \Big ) \partial_i G \Big ) \Big\|_{L_{\xi_3}^\infty ([-e^{(\chi-\delta)\tau}, e^{(\chi-\delta)\tau}]; L^2 (m)^3)} =0.
\end{split}
\end{align}
Hence, collecting \eqref{proof.sec.sec.2} and \eqref{proof.sec.sec.6}, we obtain 
\begin{align}\label{proof.sec.sec.7}
\lim_{\tau \rightarrow \infty} \sum_{|\beta|\leq 1}  \Big\| \partial_\xi^\beta \Big ( e^{\frac{\tau}{2}} \mathcal{W} ( \tau) -\sum_{i=1,2} \big ( \lambda_i + d_i \big ) \partial_i G \Big )  \Big\|_{L^\infty_{\xi_3}([-e^{(\chi-\delta)\tau},e^{(\chi-\delta)\tau}]; L^2(m)^3)} =0 
\end{align}
for any $\delta>0$, where the coefficients $d_i\in \R$ satisfies $|d_i|\leq C \| W_{0,3d} \|_{\mathbb{X}(m)}\ll 1$ with $C=C(\alpha,m,\mu,W_{0,2d})$. In particular, the last statement of Theorem \ref{theo.nlstab} holds. The proof is complete.

\section*{Acknowledgement}

The authors are grateful to the two referees for their remarks, which helped to substantially improve the manuscript.

The first author is partially supported by JSPS Program for Advancing Strategic International Networks
to Accelerate the Circulation of Talented Researchers, 'Development of Concentrated Mathematical Center Linking to Wisdom of  the Next Generation', which is organized by Mathematical Institute of Tohoku University.
The second author is partially supported by JSPS grant 25707005.
The third author acknowledges financial support from the French Agence Nationale de la Recherche under grant ANR-16-CE40-0027-01, as well as from the IDEX of the University of Bordeaux for the BOLIDE project.

\small
\bibliographystyle{abbrv}
\bibliography{lerayuloc}

\end{document}